
\def\LaTeX{%
  \let\Begin\begin
  \let\End\end
  \let\salta\relax
  \let\finqui\relax
  \let\futuro\relax}

\def\UK{\def\our{our}\let\sz s}
\def\USA{\def\our{or}\let\sz z}

\UK



\LaTeX

\USA


\salta

\documentclass[twoside,12pt]{article}
\setlength{\textheight}{24cm}
\setlength{\textwidth}{16cm}
\setlength{\oddsidemargin}{2mm}
\setlength{\evensidemargin}{2mm}
\setlength{\topmargin}{-15mm}
\parskip2mm


\usepackage[usenames,dvipsnames]{color}
\usepackage{amsmath}
\usepackage{amsthm}
\usepackage{amssymb}
\usepackage[mathcal]{euscript}
\usepackage{cite}
%
%
%


\definecolor{viola}{rgb}{0.3,0,0.7}
\definecolor{ciclamino}{rgb}{0.5,0,0.5}
\definecolor{rosso}{rgb}{0.85,0,0}

\def\juerg #1{{\color{Green}#1}}
\def\elvis #1{{\color{blue}#1}}
\def\pier #1{{\color{rosso}#1}}

\def\elvis #1{#1}
\def\pier #1{#1}
\def\juerg #1{#1}




\bibliographystyle{plain}


%

\finqui

\def\Beq{\Begin{equation}}
\def\Eeq{\End{equation}}
\def\Bsist{\Begin{eqnarray}}
\def\Esist{\End{eqnarray}}

\def\Bthm{\Begin{theorem}}
\def\Ethm{\End{theorem}}
\def\Blem{\Begin{lemma}}
\def\Elem{\End{lemma}}

\def\Bcor{\Begin{corollary}}
\def\Ecor{\End{corollary}}
\def\Brem{\Begin{remark}\rm}
\def\Erem{\End{remark}}

\def\Bcenter{\Begin{center}}
\def\Ecenter{\End{center}}
\let\non\nonumber




\def\step #1 \par{\medskip\noindent{\bf #1.}\quad}
\def\jstep #1: \par {\vspace{2mm}\noindent\underline{\sc #1 :}\par\nobreak\vspace{1mm}\noindent}


\def\aand{\quad\hbox{and}\quad}

\def\lhs{left-hand side}
\def\rhs{right-hand side}



\def\multibold #1{\def\arg{#1}%
  \ifx\arg\pto \let\next\relax
  \else
  \def\next{\expandafter
    \def\csname #1#1#1\endcsname{{\bf #1}}%
    \multibold}%
  \fi \next}

\def\pto{.}

\def\multical #1{\def\arg{#1}%
  \ifx\arg\pto \let\next\relax
  \else
  \def\next{\expandafter
    \def\csname cal#1\endcsname{{\cal #1}}%
    \multical}%
  \fi \next}


\def\multimathop #1 {\def\arg{#1}%
  \ifx\arg\pto \let\next\relax
  \else
  \def\next{\expandafter
    \def\csname #1\endcsname{\mathop{\rm #1}\nolimits}%
    \multimathop}%
  \fi \next}

\multibold
qwertyuiopasdfghjklzxcvbnmQWERTYUIOPASDFGHJKLZXCVBNM.

\multical
QWERTYUIOPASDFGHJKLZXCVBNM.

\multimathop
diag dist div dom mean meas sign supp .

\def\Span{\mathop{\rm span}\nolimits}


\def\Accorpa #1#2 #3 {\gdef #1{\eqref{#2}--\eqref{#3}}%
  \wlog{}\wlog{\string #1 -> #2 - #3}\wlog{}}


\def\separa{\noalign{\allowbreak}}

\def\supess{\mathop{\rm sup\,ess}}

\def\graffe #1{\mathopen\{#1\mathclose\}}

\def\<#1>{\mathopen\langle #1\mathclose\rangle}
\def\norma #1{\mathopen \| #1\mathclose \|}

\def\ioT {\int_0^T}
\def\intQt{\int_{Q_t}}
\def\intQ{\int_Q}
\def\iO{\int_\Omega}
\def\iG{\int_\Gamma}
\def\intS{\int_\Sigma}
\def\intSt{\int_{\Sigma_t}}

\def\bintQt{\int_{Q^t}}
\def\bintSt{\int_{\Sigma^t}}

\def\dt{\partial_t}
\def\dn{\partial_\nu}
\def\SALN{{\cal S}_{\alpha_n}^2}

\def\checkmmode #1{\relax\ifmmode\hbox{#1}\else{#1}\fi}
\def\aeO{\checkmmode{a.e.\ in~$\Omega$}}
\def\aeQ{\checkmmode{a.e.\ in~$Q$}}
\def\aeG{\checkmmode{a.e.\ on~$\Gamma$}}

\def\aet{\checkmmode{a.e.\ in~$(0,T)$}}


\def\erre{{\mathbb{R}}}

\def\nz{{\mathbb{N}}}

\def\enne{{\mathbb{N}}}




\def\genspazio #1#2#3#4#5{#1^{#2}(#5,#4;#3)}
\def\spazio #1#2#3{\genspazio {#1}{#2}{#3}T0}

\def\L {\spazio L}
\def\H {\spazio H}



\def\Lx #1{L^{#1}(\Omega)}
\def\Hx #1{H^{#1}(\Omega)}

\def\LxG #1{L^{#1}(\Gamma)}
\def\HxG #1{H^{#1}(\Gamma)}

\def\Ldue{\Lx 2}

\def\Huno{\Hx 1}
\def\Hdue{\Hx 2}

\def\HunoG{\HxG 1}
\def\HdueG{\HxG 2}

\def\LdueG{\LxG 2}



\let\theta\vartheta
\let\eps\varepsilon
\let\phi\varphi
\let\vp\varphi

\let\TeXchi\chi                         
\newbox\chibox
\setbox0 \hbox{\mathsurround0pt $\TeXchi$}
\setbox\chibox \hbox{\raise\dp0 \box 0 }
\def\chi{\copy\chibox}


\def\QED{\hfill $\square$}


\def\CO{C_\Omega}

\def\suG{_{|\Gamma}}

\def\VG{V_\Gamma}
\def\HG{H_\Gamma}
\def\WG{W_\Gamma}
\def\nablaG{\nabla_\Gamma}
\def\DeltaG{\Delta_\Gamma}
\def\muG{\mu_\Gamma}
\def\rhoG{\rho_\Gamma}
\def\tauO{\tau_\Omega}
\def\tauG{\tau_\Gamma}
\def\fG{f_\Gamma}
\def\vG{v_\Gamma}
\def\wG{w_\Gamma}

\def\gG{g_\Gamma}

\def\xiG{\xi_\Gamma}

\def\Mu{(\mu,\muG)}
\def\Rho{(\rho,\rhoG)}

\def\Xi{(\xi,\xiG)}

\def\pG{p_\Gamma}
\def\qG{q_\Gamma}

\def\gstar{g^*}

\def\rhoz{\rho_0}

\def\ual{u^\alpha}
\def\mual{\mu^\alpha}
\def\rhoal{\rho^\alpha}
\def\muGal{\mu_\Gamma^\alpha}
\def\rhoGal{\rho_\Gamma^\alpha}

\def\Mual{(\mual,\muGal)}
\def\Rhoal{(\rhoal,\rhoGal)}

\def\pal{p^\alpha}
\def\pGal{p_\Gamma^\alpha}
\def\qal{q^\alpha}
\def\qGal{q_\Gamma^\alpha}

\def\palep{p^{\alpha,\eps}}
\def\pGalep{p_\Gamma^{\alpha,\eps}}
\def\qalep{q^{\alpha,\eps}}
\def\qGalep{q_\Gamma^{\alpha,\eps}}

\def\Pi{\hat\pi}

\def\piG{\pi_\Gamma}

\def\un{u^{\alpha_n}}
\def\mun{\mu^{\alpha_n}}
\def\muGn{\mu_\Gamma^{\alpha_ n}}
\def\rhon{\rho^{\alpha_n}}
\def\rhoGn{\rho_\Gamma^{\alpha_n}}
\def\Rhon {(\rhon,\rhoGn)}
\def\Mun {(\mun,\muGn)}

\def\calVz{\calV_0}
\def\calHz{\calH_0}
\def\calVsz{\calV_{*0}}
\def\calVzp{\calV_0^{\,*}}
\def\calVp{\calV^{\,*}}

\def\normaHH #1{\norma{#1}_{\calH}}
\def\normaVV #1{\norma{#1}_{\calV}}

\let\hat\widehat

\def\uad{{\cal U}_{\rm ad}}


\Begin{document}


%
\title{Optimal velocity control of a convective Cahn--Hilliard \\[3mm] system with 
double obstacles  
  and dynamic  boundary \\[3mm] conditions: a `deep quench' approach\footnote{PC and GG gratefully acknowledge finantial support from the MIUR-PRIN Grant 2015PA5MP7 
``Calculus of Variations'', the GNAMPA (Gruppo Nazionale per l'Analisi
Matematica, la Probabilit\`{a} e loro Applicazioni) of INDAM (Istituto Nazionale
di Alta Matematica) and the IMATI -- C.N.R. Pavia.}}
\author{}
\date{}
\maketitle
\Bcenter
\vskip-2cm
{\large\sc Pierluigi Colli$^{(1)}$}\\
{\normalsize e-mail: {\tt pierluigi.colli@unipv.it}}\\[.25cm]
{\large\sc Gianni Gilardi$^{(1)}$}\\
{\normalsize e-mail: {\tt gianni.gilardi@unipv.it}}\\[.25cm]
{\large\sc J\"urgen Sprekels$^{(2)}$}\\
{\normalsize e-mail: {\tt sprekels@wias-berlin.de}}\\[.45cm]
$^{(1)}$
{\small Dipartimento di Matematica ``F. Casorati'', Universit\`a di Pavia}\\
{\small and Research Associate at the IMATI -- C.N.R. Pavia}\\
{\small via Ferrata 5, 27100 Pavia, Italy}\\[.2cm]
$^{(2)}$
{\small Department of Mathematics}\\
{\small Humboldt-Universit\"at zu Berlin}\\
{\small Unter den Linden 6, 10099 Berlin, Germany}\\[2mm]
{\small and}\\[2mm]
{\small Weierstrass Institute for Applied Analysis and Stochastics}\\
{\small Mohrenstrasse 39, 10117 Berlin, Germany}
\Ecenter
\Begin{abstract}
\noindent In this paper, we investigate a distributed optimal control problem for a 
convective viscous Cahn--Hilliard system with dynamic boundary conditions.
Such systems govern phase separation processes between two phases taking place in an incompressible fluid in a container and, at the same time, on the container boundary.  
The cost functional is of standard tracking type, while the control is exerted
by the velocity of the fluid in the bulk. In this way, the coupling between the state 
(given by the associated order parameter and chemical potential) 
and control variables in the governing system of nonlinear partial differential equations is bilinear, 
which presents a difficulty for the analysis. 
In contrast to the previous paper 
\pier{{\em Optimal velocity control of a viscous CahnÐHilliard system with convection and dynamic boundary conditions} by the same authors,} 
the bulk and surface free energies are of double obstacle type, which renders the state constraint
nondifferentiable. It is well known that for such cases standard constraint qualifications
are not satisfied so that standard methods do not apply to yield the existence of 
Lagrange multipliers. In this paper, we overcome this difficulty by taking advantage
of results established in the quoted paper for logarithmic nonlinearities, using a so-called
`deep quench approximation'. We derive results concerning the existence of optimal
controls and the first-order necessary optimality conditions in terms of a variational
inequality and the associated adjoint system. 
     

\End{abstract}
\salta
\pagestyle{myheadings}
\newcommand\testopari{\sc Colli \ --- \ Gilardi \ --- \ Sprekels}
\newcommand\testodispari{\sc Velocity control of convective Cahn--Hilliard system}
\markboth{\testopari}{\testodispari}
\finqui
%

\section{Introduction}
\label{INTRO}
\setcounter{equation}{0}

Let $\Omega\subset\erre^3$ denote some open, bounded and connected set having a smooth boundary $\Gamma$
and unit outward normal $\,\nu$.
We denote by $\dn$, $\nablaG$, $\Delta_\Gamma$ the outward normal derivative, the tangential gradient,
and the Laplace--Beltrami operator on $\Gamma$, in this order. Moreover, we fix some final time $T>0$ and
introduce for every $t\in (0,T]$ the sets $Q_t:=\Omega\times (0,t)$ and $\Sigma_t:=\Gamma\times (0,t)$,
where we put, for the sake of brevity, $Q:=Q_T$ and $\Sigma:=\Sigma_T$. We then consider the following optimal control
problem:

\vspace{3mm}\noindent
(${\cal P}_0$) \quad Minimize the cost functional
\begin{align} 
\label{cost} 
{\cal J}((\rho, \rhoG),u)
   :=&\,  \frac{\beta_1}2 \intQ |\rho-\widehat \rho_Q|^2
  + \frac{\beta_2}2 \intS |\rho-\widehat \rho_\Sigma|^2
  \non
  \\
  & 
  \,+ \frac{\beta_3}2 \iO |\rho(T)-\widehat\rho_\Omega|^2
  + \frac{\beta_4}2 \iG |\rhoG(T)-\widehat\rho_\Gamma|^2
  + \frac{\beta_5}2 \intQ |u|^2\,,
\end{align} 
subject to  the state system  
\begin{align}
\label{ss1}
&\dt\rho+\nabla\rho\cdot u-\Delta\mu=0 \quad\mbox{in }\,Q\,,\\[1mm]
\label{ss2}
&\tauO\dt\rho-\Delta\rho+\xi+\pi(\rho)=\mu \quad\mbox{in }\,Q\,,\\[1mm]
\label{ss3}
&\xi\in\partial I_{[-1,1]} \juerg{(\rho)}\quad\mbox{in }\,Q\,,\\[1mm]
\label{ss4}
&\dt\rhoG+\dn\mu-\DeltaG\muG=0 \quad\mbox{and}\quad \mu_{|\Sigma}=\muG\quad\mbox{on }\,\Sigma\,,\\[1mm]
\label{ss5}
&\tauG\dt\rhoG+\dn\rho-\DeltaG\rhoG+\xiG+\piG(\rhoG)=\muG \quad\mbox{and}\quad \rho_{|\Sigma}=\rhoG
\quad\mbox{on }\,\Sigma\,,\\[1mm]
\label{ss6}
&\xiG\in \partial I_{[-1,1]} \juerg{(\rhoG)}\quad\mbox{on }\,\Sigma\,,\\[1mm]
\label{ss7}
&\rho(0)=\rho_0\quad\mbox{in }\,\Omega,\quad \rhoG(0)=\rho_{0|\Gamma}\quad\mbox{on }\,\Gamma\,, 
\end{align}
and to the control constraint
\begin{equation}
\label{concon}
u\in\uad\,.
\end{equation}
Here, the constants $\beta_i$, $1\le i\le 5$, are nonnegative but not all zero, and 
 $\widehat\rho_Q$, $\widehat \rho_\Sigma$, $\widehat\rho_\Omega$, 
$\widehat\rho_\Gamma$, are given target functions.
Furthermore, $\pi$, $\piG$ denote smooth functions,
{while $I_{[-1,1]}$ is the indicator function of the interval~$[-1,1]$}.
Moreover, $\uad$ is a suitable bounded, closed and convex subset of the control space
\begin{equation}
\label{conspace}
\calX:= L^2(0,T;\widetilde U)\cap (L^\infty(Q))^3\cap (H^1(0,T;L^3(\Omega)))^3\,,
\end{equation}
where
\begin{equation}
\label{deftilU}
\widetilde U\,:=\,\left\{u\in (L^2(\Omega))^3: \,\mbox{${\rm div}\,u=0$ \,a.e. in\, $\Omega$ \,and\, 
$u\cdot\nu=0$ \,a.e. on \,$\Gamma$}\right\}. 
\end{equation}
\noindent
The regularity condition $\,u\in (H^1(0,T;L^3(\Omega)))^3\,$ for the admissible controls seems to be unusual 
at a first glance. However, in view of the bilinear coupling between control and state,
it turns out (cf.~\cite{CGS13}) that, among other constraints, this is exactly the kind of regularity 
that guarantees the existence of a unique solution to the state
system having sufficient regularity properties. 

 We note that the state system \eqref{ss1}--\eqref{ss7} can be seen as a phase field model 
 for a  phase separation process taking place in an incompressible fluid in the container $\Omega$ and
 on the container boundary $\Gamma$. In this connection, the variables $\Mu$ and $\Rho$ stand
 for the chemical potential and the order parameter (usually the density of one of the involved
 phases, normalized in such a way as to attain its values in the interval [-1,1]) 
 of the phase separation process in the bulk and on the surface, respectively. It is worth
 noting that the total mass of the order parameter is conserved during the separation process;
 indeed, integrating \eqref{ss1} for fixed $t\in (0,T]$ over $\Omega$, and using the fact that
 $\,u\in\calX\,$, as well as \eqref{ss4}, we readily find that
 \begin{equation}
 \label{conserve}
 \dt\Big(\iO\rho(t)+\iG\rhoG(t)\Big)=0\,.
 \end{equation} 
 We also note that the densities of the local free bulk energy $f+I_{[-1,1]}$ and the
local free surface energy $\fG+I_{[-1,1]}$ are typically of double obstacle type.

\pier{In the mathematical literature numerous contributions are dedicated to the questions of well-posedness and asymptotic behavior
for various types of Cahn--Hilliard systems: viscous or nonviscous, local or nonlocal, 
with zero Neumann boundary conditions or dynamic boundary conditions. We omit to (try to) 
quote a number of contributions since they are too many and we would surely miss some of the 
important ones. However, let us point out that there are
still a few papers dealing with the related optimal control problems: among them, we
refer to \cite{CGRS1,CGS2,CGS4,DZ,hw,WaNa,Zh,ZW} for the case of Dirichlet or zero 
Neumann boundary conditions and to \cite{CFGS1,CFGS2,CGS1,CGS3, CGS3bis,CS10, FY} for the case of dynamic boundary conditions.}

\pier{A recent investigation for convective Cahn--Hilliard systems produced 
the results rigorously proved in\cite{ZL1} for the one-dimensional and 
in\cite{ZL2} for the two-dimensional case. The papers\cite{FRS,TM} are devoted to
the distributed optimal control of a two-dimensional 
Cahn--Hilliard/Navier--Stokes system. Let us also mention the contributions \cite{HHKW,HKW,HW1,HW2}, which deal with the optimal control of the 
Cahn--Hilliard/Navier--Stokes system in 3D, but for some time-discretized version.}

\pier{A key feature of this paper is the use of the fluid velocity as the control variable in the convective Cahn--Hilliard system. From a practical point of view, this control process can be realized by placing either a mechanical stirring device or 
an ultrasound emitter into the container. In the case of electrically conducting fluids
like molten metals, a remarkable option is the possibility of using magnetic fields 
(cf.~\cite{Kudla} for applications of this kind).} 
To the authors' best knowledge, the only existing mathematical contributions, 
in which the fluid velocity
is used as the control in a convective Cahn--Hilliard system in three dimensions of space, 
are the recent contributions \cite{RS} and \cite{CGS14}. While in \cite{RS} a nonlocal convective 
Cahn--Hilliard system with a possibly degenerating mobility and zero
Neumann boundary conditions was studied, we considered in \cite{CGS14} a viscous local Cahn--Hilliard 
system with constant mobility (normalized to unity) and the more difficult dynamic
boundary conditions (see also \cite{CGS13} and \cite{CF2} for related results). However, in \cite{CGS14} only differentiable nonlinearities were admitted.

In this contribution, we investigate the much more challenging nondifferentiable double obstacle 
case when $\xi,\xi_\Gamma$ satisfy the inclusions (\ref{ss3}), \eqref{ss6}, and we assume dynamic 
boundary conditions. Moreover, we consider the 
spatially three-dimensional case.  Our approach is guided by a strategy that was introduced  
by two of the present authors and M. H. Farshbaf-Shaker in \cite{CFS}: we aim to derive first-order necessary
optimality conditions for the double obstacle case by 
performing a so-called `deep quench limit' in a family of optimal control problems
with differentiable logarithmic nonlinearities that  was treated in \cite{CGS14},
and for which the corresponding state systems were analyzed in \cite{CGS13}.
The general idea is briefly explained as follows: we replace the inclusions (\ref{ss3}) and
\eqref{ss6} by the identities
\begin{equation}
\label{ss4neu}
\xi=\vp(\alpha)\,h'(\rho) \pier{\quad\hbox{in } Q}, \quad \xi_\Gamma=\vp(\alpha)\,h'(\rhoG)\pier{\quad\hbox{on } \Sigma},
\end{equation}
where $h$ is defined by
\begin{equation}
\label{defh}
h(\rho):=\left\{\begin{array}{ll}
(1-\rho)\,\ln(1-\rho)+(1+\rho)\,\ln(1+\rho)&\mbox{if } \,\rho\in (-1,1)\\[1mm]
2\,\ln(2)&\mbox{if }\,\rho\in\{-1,1\}
\end{array}\right.,
\end{equation}
and where 
\begin{equation}
\label{phiat0}
\mbox{ $\varphi\in C(0,1]$ is positive on $(0,1]$ and satisfies }\,
\lim_{\alpha\searrow 0}\,\varphi(\alpha)=0.
\end{equation}
We remark that we can simply choose $\,\varphi(\alpha)=\alpha^p\,$ for some $\,p>0$.
Now observe that \pier{$h(y)\geq 0$ for all $y\in [-1,1]$,}
$h'(y)=\ln\left(\frac{1+y}{1-y}\right)$ \,and\, $h''(y)=\frac 2 {1-y^2}>0$\, for 
$y\in (-1,1)$. Hence, in particular, we have
\begin{align}
\label{hphi}
&\pier{\lim_{\alpha\searrow 0}\,\varphi(\alpha)\,h(y)=0 \quad\forall\, y\in [-1,1],}
\quad 
\lim_{\alpha\searrow 0}\,\varphi(\alpha)\,h'(y)=0  \pier{\quad\forall\, y\in (-1,1),} \nonumber\\[2mm]
&\lim_{\alpha\searrow 0}\Bigl(\varphi(\alpha)\,\lim_{y\searrow -1}h'(y)\Bigr)\,=\,-\infty,
\quad \lim_{\alpha\searrow 0}\Bigl(\varphi(\alpha)\,\lim_{y\nearrow +1}h'(y)\Bigr)\,=\,+\infty\,.
\end{align}
We thus may regard the graph 
$\,\varphi(\alpha)\,h'\,$ as an approximation to the graph of the subdifferential
$\partial I_{[-1,1]}$. 

\vspace{2mm}\quad
Now, for any $\alpha>0$, the optimal control problem (later to be denoted by $({\cal P}_\alpha)$), 
which results if in $({\cal P}_0)$ the relations (\ref{ss3}), \eqref{ss6} are replaced by 
(\ref{ss4neu}), is of the type for which
in \cite{CGS14} the existence of optimal controls $u^\alpha\in\uad$ as well as first-order 
necessary optimality conditions have been derived. Proving a priori estimates (uniform in $\alpha>0$), and 
employing compactness and monotonicity arguments, we will be able to show the following existence and 
approximation result: whenever $\,\{u^{\alpha_n}\}\subset\uad$ is a sequence of optimal controls for 
$({\cal P}_{\alpha_n})$, where $\alpha_n\searrow 0$ as $n\to\infty$, then there exist
a subsequence of $\{\alpha_n\}$, which is again indexed by $n$, and an optimal control 
$\bar u\in\uad$ of
$({\cal P}_0)$ such that
\begin{equation}
\label{eq:1.17}
u^{\alpha_n}\to\bar u \quad\mbox{weakly-star in ${\cal X}$ as }\,
n\to\infty\,.
\end{equation}
In other words, optimal controls for $({\cal P}_\alpha)$ are for small $\alpha>0$ likely to be `close' to 
optimal controls for $({\cal P}_0)$. It is natural to ask if the reverse holds, i.\,e., whether every optimal control for
 $({\cal P}_0)$ can be approximated by a sequence $\,\{u^{\alpha_n}\}\,$ of optimal controls
for $({\cal P}_{\alpha_n})$, for some sequence $\alpha_n\searrow 0$. 

\vspace{2mm}\quad
Unfortunately, we will not be able to prove such a `global' result that applies to all optimal controls for
(${\cal P}_0$). However,  a `local' result can be established. To this end, let $\bar u\in\uad$ be any optimal control
for $({\cal P}_0)$. We introduce the `adapted' cost functional
\begin{equation}
\label{adcost}
\widetilde{\cal J}((\rho,\rhoG),u) \,:=\,{\cal J}((\rho,\rhoG),u)\,+\,\frac 1 2\,
\|u-\bar u\|^2_{(L^2(Q))^3}
\end{equation}
and consider for every $\alpha\in (0,1]$ the {\em adapted control problem} of minimizing 
$\,\widetilde{\cal J}\,$ subject to $u\in\uad$ and to the constraint that $((\mu,\muG),(\rho,\rhoG))$ 
solves the approximating system (\ref{ss1}), (\ref{ss2}), \eqref{ss4}, (\ref{ss5}), \eqref{ss7}, 
(\ref{ss4neu}). It will then turn out that the following is true: 

\vspace{2mm}\noindent
(i) \,There are some sequence $\,\alpha_n\searrow 0\,$ and minimizers 
$\,{\bar u^{\alpha_n}}\in\uad$ of the adapted control problem 
associated with $\alpha_n$, $n\in\nz$,
such that
\begin{equation}
\label{eq:1.19}
{\bar u^{\alpha_n}}\to\bar u\quad\mbox{strongly in $(L^2(Q))^3$
as }\,n\to \infty.
\end{equation}
(ii) It is possible to pass to the limit as $\alpha\searrow 0$ in the first-order necessary
optimality conditions corresponding to the adapted control problems associated with $\alpha\in 
(0,1]$ in order to derive first-order necessary optimality conditions for problem $({\cal P}_0)$.

\vspace{2mm}\quad
The paper is organized as follows: in Section~2, we give a precise statement of the problem
under investigation, and we derive some results concerning the state system 
(\ref{ss1})--(\ref{ss7}) and 
its $\alpha\,$--\,approximation which is obtained if in $({\cal P}_0)$ the relations
 (\ref{ss3}) and \eqref{ss6} are replaced by the relations (\ref{ss4neu}).
In Section~3, we then prove the existence of optimal controls and the approximation 
result formulated above in
(i). The final Section~4 is devoted to the derivation of the first-order necessary 
optimality conditions, where the  strategy outlined in (ii) is employed. 

\vspace{2mm}
During the course of this analysis, we will make 
repeated use of H\"older's inequality, of the elementary Young's inequality
\begin{equation}
\label{Young}
a\,b\,\le\,\gamma |a|^2\,+\,\frac 1{4\gamma}\,|b|^2\quad\forall\,a,b\in\erre, \quad\forall\,\gamma>0,
\end{equation}
as well as the continuity of the embeddings $H^1(\Omega)\subset L^p(\Omega)$ for $1\le p\le 6$ and 
$\Hdue\subset C^0(\overline\Omega)$. Notice that the latter embedding is also compact, while this holds true
for the former embeddings only if $p<6$. We will also use the denotations
\begin{align}
\label{defQt}
&Q^t:=\Omega\times (t,T),\quad \Sigma^t:=\Gamma\times (t,T),\quad\mbox{for }\,0\le t<T.
\end{align}
Moreover, throughout the paper, \pier{for a Banach space $\,X\,$ we denote by $\,X^*\,$ its dual space.  Let $\,\|\,\cdot\,\|_X\,$ stand for the norm in the space $X$ or in a
power of it.} The only exemption from this rule \pier{is for} the norms of the
$\,L^p\,$ spaces and of their powers, which we often denote by $\|\,\cdot\,\|_p$, for
$\,1\le p\le +\infty$. By $\,\langle v,w\rangle_X\,$ we will always denote
the dual pairing between elements $\,v\in X^*\,$ and $\,w\in X$. Finally, we recall 
some well-known estimates from trace theory and from the theory of
elliptic equations. Namely, there is some constant $C_\Omega>0$, which depends only on
$\Omega$, such that, for every $v$ and $\vG$ for which the \rhs s are meaningful,
\begin{align}
\label{CO1}
\|v\|_{H^{3/2}(\Omega)}&\,\le\,\CO \left(\|v_{|\Gamma}\|
_{H^1(\Gamma)}\,+\,\|\Delta v\|_{\pier{L^2(\Omega)}}\right),\\[1mm]
\label{CO2}
\|\dn v\|_{L^2(\Gamma)}&\,\le\,\CO \left(\|v\|_{H^{3/2}(\Omega)}
\,+\,\|\Delta v\|_{\pier{L^2(\Omega)}}\right),\\[1mm]
\label{CO3}
\|v_\Gamma\|_{H^2(\Gamma)}&\,\le\,\CO \left(\|v_\Gamma\|
_{H^1(\Gamma)}\,+\,\|\DeltaG\vG\|_{L^2(\Gamma)}\right),\\[1mm]
\label{CO4}
\|v\|_{H^2(\Omega)}&\,\le\,\CO \left(\|v_{|\Gamma}\|
_{H^{3/2}(\Gamma)}\,+\,\|\Delta v\|_{\pier{L^2(\Omega)}}\right).
\end{align}


\section{General setting and the state system}
\label{STATE}
\setcounter{equation}{0}

In this section, we introduce the general setting of our control 
problem and state some results on the state system \eqref{ss1}--\eqref{ss7}. To begin with, we
recall the definition~\eqref{conspace} of~$\calX$ and introduce the spaces
\begin{align}
  & H := \Ldue \,, \quad  
  V := \Huno\,,   \quad
  W := \Hdue,
  \label{defspaziO}
  \\
  & \HG := \LdueG \,, \quad 
  \VG := \HunoG \,,
  \quad
  \WG := \HdueG,
  \label{defspaziG}
  \\
  & \calH := H \times \HG \,, \quad
  \calV := \graffe{(v,\vG) \in V \times \VG : \ \vG = v\suG}\,,
  \quad 
  \calW := \bigl( W \times \WG \bigr) \cap \calV \,.
  \label{defspaziprod}
\end{align}
In the following, we will often work in the framework of the Hilbert triplet
$(\calV,\calH,\calVp)$.
Thus, we have
$$\<(g,\gG),(v,\vG)>_{\calV}=\iO gv+\iG\gG\vG \quad\mbox{
for every $(g,\gG)\in\calH$ and $(v,\vG)\in\calV$}.$$
Next, denote by $\,(1,1)\in\calV$\, the pair whose component functions equal unity in $\Omega$ and on~$\Gamma$,
respectively,  and by $\,|\Omega|\,$
and $\,|\Gamma|\,$ the volume of $\Omega$ and the area of $\Gamma$, respectively. 
We then define 
the generalized mean value of a functional $\,g^*\in\calVp$ by
\Beq
  \mean \gstar := \frac { \< \gstar,(1,1) >_{\calV} } { |\Omega| + |\Gamma| },
  \label{genmean}
\Eeq
which, if $\gstar=(v,\vG)\in\calH$, becomes
\Beq
  \mean\,(v,\vG) = \frac { \iO v + \iG \vG } { |\Omega| + |\Gamma| }.
  \label{usualmean}
\Eeq
Observe that the function
\Beq
  \calV \ni (v,\vG) \mapsto \int_\Omega\pier{|\nabla v|}^2 + \iG|\nablaG\vG|^2 +   |\mean(v,\vG)|^2
  \non
\Eeq
yields the square of a Hilbert norm on $\calV$ that is equivalent to the natural one,
i.e., we have, \pier{for} some $\CO>0$ which depends only on $\Omega$,
\Beq
  \normaVV{(v,\vG)}^2 \leq \CO \Big(\iO |\nabla v|^2\, + \iG|\nablaG\vG|^2\, +\, |\mean(v,\vG)|^2\Big)
  \quad \forall\, (v,\vG)\in\calV\,.
  \label{normaVVequiv}
\Eeq
Next, we set
\Beq
  \calVsz := \graffe{ \gstar\in\calVp : \ \mean\gstar = 0 } , \quad
  \calHz := \calH \cap \calVsz
  \aand
  \calVz := \calV \cap \calVsz .
  \label{defcalVz}
\Eeq
\Accorpa\Defspazi defspaziO defcalVz
Notice the difference between $\calVsz$ and the dual space~$\calVzp=(\calVz)^*$.
At this point, it is clear that the function
\Beq
  \calVz \ni (v,\vG) \mapsto \norma{(v,\vG)}_{\calVz}^2
  := \iO|\nabla v|^2+\iG|\nablaG\vG|^2
  \label{normaVz}
\Eeq
is the square of a Hilbert norm on $\calVz$ which is equivalent to the usual one.
This has the consequence (see\pier{\cite[Sect.~2]{CGS13}})
that, for every $\gstar\in\calVsz$, there exists a unique pair 
$\Xi\in\calVz$ such that
\Beq
  \iO \nabla\xi \cdot \nabla v + \iG \nablaG\xiG \cdot \nablaG\vG
  = \< \gstar , (v,\vG) >_{\calV}
  \quad \hbox{for every $(v,\vG)\in\calV$}.
  \qquad
  \label{perdefN}
\Eeq
This allows us to define the operator $\calN:\calVsz\to\calVz$ as follows:
\Beq
  \hbox{For $\gstar\in\calVsz$, \ $\calN\gstar$ is the unique pair $\Xi\in\calVz$ satisfying \eqref{perdefN}}.
  \label{defN}
\Eeq
We notice that $\calN$ is linear, symmetric, and bijective.
Therefore, if we set
\Beq
  \norma\gstar_* := \norma{\calN\gstar}_{\calVz},
  \quad \hbox{for $\gstar\in\calVsz$},
  \label{normastar}
\Eeq
then we obtain a Hilbert norm on $\calVsz$ which 
turns out to be equivalent to the norm induced by the norm of~$\calVp$.
For future use, we collect some properties of~$\calN$.
By just applying the definition, we readily see that
\begin{align}
  & \< \gstar , \calN\gstar >_{\calV}
  = \norma\gstar_*^2
  \quad\, \hbox{if $\gstar \in \calVsz$}\,,
  \label{propNa}
  \\[1mm]
  & \iO \nabla w \cdot \nabla\xi
  + \iG \nablaG\wG \cdot \nablaG\xiG
  = \normaHH{(w,\wG)}^2
  \quad\, \hbox{if $(w,\wG) \in \calVz$ and $(\xi,\xiG) = \calN(w,\wG)$} \,.
  \label{propNb}
\end{align}
Moreover, owing to the symmetry of~$\calN$ (where, here and in the following, $\calN$~is also applied to $\calVsz$-valued functions in the obvious way), we have, for a.e. $t\in (0,T)$, that
\begin{align}
  & \langle \dt\gstar(t) , \calN\gstar(t)\rangle_{\calV}
  =  \frac 12 \, \frac d{dt} \, \|\gstar(t)\|_*^2\,,\quad\, \hbox{if $\gstar \in \H1\calVsz$},
  \label{propNdta}
  \\[2mm]
  & \iO \nabla w(t) \cdot \nabla\xi(t)
  + \iG \nablaG\wG(t) \cdot \nablaG\xiG(t)
  = \frac 12 \, \frac d{dt} \, \normaHH{(w(t),\wG(t))}^2\,,
  \non
  \\[1mm]
  &  \hbox{\,\,if $(w,\wG)\in\L2\calV$, \ $\dt(w,\wG) \in \L2\calVsz$\,    and\, 
  $\Xi = \calN(\dt(w,\wG))$} \,.
  \label{propNdtb}
\end{align}
\Accorpa\PropN propNa propNdtb

\vspace{5mm}
We now turn our interest to the state system \eqref{ss1}--\eqref{ss7}, observing that with the 
above notations its weak form reads as follows:
we look for functions $(\Mu,\Rho)$ such that $\,\mu_{|\Sigma}=\muG\,$ and \,$\rho_{|\Sigma}
=\rhoG\,$ as well as 
\begin{align}
\label{prima}
  & \iO \dt\rho \, v
  + \iG \dt\rhoG \, \vG
  - \iO \rho u \cdot \nabla v
  + \iO \nabla\mu \cdot \nabla v
  + \iG \nablaG\muG \cdot \nablaG\vG
  = 0
  \non
  \\[1mm]
  & \quad \hbox{\aet\ and for every $(v,\vG)\in\calV$},
   \\[2mm]
  \label{seconda}
  & \tauO \iO \dt\rho \, v
  + \tauG \iG \dt\rhoG \, \vG
  + \iO \nabla\rho \cdot \nabla v
  + \iG \nablaG\rhoG \cdot \nablaG\vG
  \non
  \\[1mm]
  & \quad {}
  + \iO (\xi+{\pi}(\rho)) v
  + \iG (\xiG+{\piG}(\rhoG)) \vG
  = \iO \mu v 
  + \iG \muG \vG
  \non
  \\[1mm]
  & \quad \hbox{\aet\ and for every $(v,\vG)\in\calV$},
     \\[2mm]
 \label{terza}    
&\xi\in\partial I_{[-1,1]}(\rho) \quad\mbox{a.e. in }\,\juerg{Q},\qquad
\xiG\in\partial I_{[-1,1]}(\rhoG) \quad\mbox{a.e. on }\,\Sigma,\\[2mm]    
  & \rho(0) = \rhoz
  \quad \aeO, \qquad \rhoG(0)=\rho_{0|\Gamma}\quad \aeG .
  \label{cauchy}
\end{align}

We make the following assumptions on the data of our problem:
\begin{description}
\item{(A1)} \quad$(\rho_0,\rho_{0|\Gamma}) \in\calW$, and we have $\,-1<\rho_0(x)<1\,$ for all $\,x\in\overline\Omega$. 
\item{(A2)} \quad $\tauO>0$ and $\tauG>0$.
\item{(A3)} \quad {$\pi,\piG\in C^2[-1,1]$}.
\item{(A4)} \quad The constants $\,\beta_i$, $1\le i\le 5$, are all nonnegative but not all
 equal to zero, and \hspace*{5mm} it holds $\,\widehat \rho_Q\in L^2(Q)$, $\,\widehat \rho_\Sigma
\in L^2(\Sigma)$, $\,\widehat\rho_\Omega\in\Ldue$, and $\,\widehat\rho_\Gamma\in L^2(\Gamma)$. 
\item{(A5)} \quad
The function $\overline U\in L^\infty(Q)$ and the constant $R_0>0$ make the 
admissible set
\begin{align}
 \uad := &\left\{ u\in{\cal X}: \,|u|\leq\overline U \ \aeQ,\, 
\|u\|_{\calX}\leq R_0 \right\}
   \label{defUad}
\end{align}
\hspace*{6mm}nonempty.
\end{description}

\vspace{3mm}
\Brem
Notice that the conditions  $\,{\rm div}\,u=0\,$ in $\Omega$, \,$u\cdot\nu=0$\, on $\Gamma$,
encoded in the definition of $\calX$, have
to be understood in the generalized sense, i.e., they are equivalent to postulating that
\begin{equation}
\label{divu0}
\iO u\cdot \nabla v=0 \quad \forall\,v\in V.
\end{equation}
We thus may infer that $\uad$ is a bounded, closed and convex subset of $\calX$. 
\Erem	

The following result is a special case of \cite[Thms. 2.3,~2.6]{CGS13}.

\Bthm
\label{WPstate}
Suppose that the assumptions {\rm (A1)--(A3)} and {\rm (A5)} hold true. 
Then the state system {\rm \eqref{ss1}--\eqref{ss7}}  has for every $u\in\uad$ at least one solution
$(\Mu,\Rho,(\xi,\xiG))$ such that
\begin{align}
\label{regmu}
&\Mu\in 
L^\infty(0,T;\calW),\\[1mm]
\label{regrho}
&\Rho\in W^{1,\infty}(0,T;\calH)\cap H^1(0,T;\calV)\cap L^\infty
(0,T;\calW),\\[1mm]
\label{regxi}
&(\xi,\xiG)\in L^\infty(0,T;\calH).
 \end{align}
Moreover, the component $\Rho$ is the same for any such solution. In addition, there is some
constant $\,K_1^*>0$, which depends only on the data of the problem, such that for any solution 
$(\Mu,\Rho,(\xi,\xiG))$ associated with some $\,u\in\uad\,$ it holds that
\begin{align}
\label{ssb1}
&\|\Mu\|_{L^\infty(0,T;\calW)}\,+\,\|\Rho\|_{W^{1,\infty}(0,T;\calH)\cap H^1(0,T;\calV)\cap L^\infty
(0,T;\calW)}\non\\[1mm]
&+\,\|(\xi,\xiG)\|_{L^\infty(0,T;\calH)}\,\le\,K_1^*\,.
\end{align}
\Ethm

It follows from Theorem~2.2 that the mapping $\uad\ni u\mapsto \calS^2_0(u):=\Rho$ is well defined. Next, we consider for $\alpha\in (0,1]$ the $\alpha-$approximating system
\begin{align}
\label{as1}
&\dt\rhoal+\nabla\rhoal\cdot u-\Delta\mual=0 \quad\mbox{a.e. in }\,Q\,,\\[1mm]
\label{as2}
&\tauO\,\dt\rhoal-\Delta\rhoal+\vp(\alpha)h'(\rhoal)+{\pi}(\rhoal)=\mual \quad\mbox{a.e. in }\,Q\,,\\[1mm]
\label{as3}
&\dt\rhoGal+\dn\mual-\DeltaG\muGal=0 
\quad\mbox{and}\quad \mu^\alpha_{|\Sigma}=\muGal\quad\mbox{a.e. on }\,\Sigma\,,\\[1mm]
&\tauG\,\dt\rhoGal+\dn\rhoal-\DeltaG\rhoGal+\vp(\alpha)h'(\rhoGal)+{\piG}(\rhoGal)=\muGal
\non
\\[1mm]
\label{as4}
& \quad{\mbox{and}\quad \rho^\alpha_{|\Sigma}=\rhoGal
\quad\mbox{a.e. on }\,\Sigma\,,}\\[1mm]
\label{as5}
&\rhoal(0)=\rho_0\quad\mbox{a.e. in }\,\Omega,\quad \rhoGal(0)=\rho_{0|\Gamma}\quad\mbox{a.e. on }\,\Gamma\,, 
\end{align}
where $\,h\,$ is given by \eqref{defh} and $\,\vp\,$ satisfies \eqref{phiat0}.
The corresponding weak formulation reads as follows: we look for functions $(\Mual,\Rhoal)$ such that
$\,\mu^\alpha_{|\Sigma}=\muGal\,$ and $\,\rho^\alpha_{|\Sigma}=\rhoGal\,$ as well as
\begin{align}
\label{was1}
  & \iO \dt\rhoal \, v
  + \iG \dt\rhoGal \, \vG
  - \iO \rhoal u \cdot \nabla v
  + \iO \nabla\mual \cdot \nabla v
  + \iG \nablaG\muGal \cdot \nablaG\vG
  = 0
  \non
  \\[1mm]
  & \quad \hbox{\aet\ and for every $(v,\vG)\in\calV$},
   \\[2mm]
  \label{was2}
  & \tauO \iO \dt\rhoal \, v
  + \tauG \iG \dt\rhoGal \, \vG
  + \iO \nabla\rhoal \cdot \nabla v
  + \iG \nablaG\rhoGal \cdot \nablaG\vG
  \non
  \\[1mm]
  & \quad {}
  + \iO (\vp(\alpha)h'(\rhoal)+{\pi}(\rhoal)) v
  + \iG (\vp(\alpha)h'(\rhoGal)+{\piG}(\rhoGal)) \vG
  = \iO \mual v 
  + \iG \muGal \vG
  \non
  \\[1mm]
  & \quad \hbox{\aet\ and for every $(v,\vG)\in\calV$},
     \\[2mm]
 \label{was3}    
  & \rhoal(0) = \rhoz
  \quad \aeO, \qquad \rhoGal(0)=\rho_{0|\Gamma}\quad \aeG .
  \end{align}
Observe that also this system has the property that the unknown representing the order parameter is a conserved
quantity: indeed, insertion of $(v,\vG)=(1,1)\in\calV$ in \eqref{was1} and integration over time yield
 that
\begin{equation}
\label{meanal}
\hat r\,:=\,{\rm mean}\,(\rhoz,\rho_{0|\Gamma})\,=\,{\rm mean}\,(\rhoal(t),\rhoGal(t)), \mbox{ $0\le t\le T$,} 
\quad\forall \,\alpha\in (0,1].
\end{equation}

We have the following result for the approximating system.
\Bthm
Suppose that the conditions {\rm (A1)--(A3)}, {\rm (A5)}, \eqref{defh} and \eqref{phiat0}
  are satisfied. Then the system 
{\rm \eqref{as1}--\eqref{as5}}   
has for every $\alpha\in (0,1]$ and for every $u\in\uad$ a unique solution $(\Mual,\Rhoal)$ satisfying {\rm \eqref{regmu}}
and {\rm \eqref{regrho}}.
Moreover, there are constants $\rho_*(\alpha),\rho^*(\alpha)\in (-1,1)$ and \,$K_2^*>0$, which depend only on the data
of the state system, such that the following
holds true: whenever $(\Mual,\Rhoal)$ is the solution to the system {\rm \eqref{as1}--\eqref{as5}} associated with some $\,\alpha\in (0,1]$\, and \,$u\in\uad$,
then we have
\begin{align}
\label{separ}
&\rho_*(\alpha)\,\le\rhoal(x,t)\,\le\rho^*(\alpha) \quad\forall\,(x,t)\in \overline Q\,,\\[2mm]
\label{wsb1}
&\|\Mual\|_{L^\infty(0,T;\calW)}\,+\,\|\Rhoal\|_{W^{1,\infty}(0,T;\calH)\cap H^1(0,T;\calV)\cap L^\infty
(0,T;\calW)}\non\\[1mm]
&+\,\|(\vp(a)h'(\rhoal),\vp(\alpha)h'(\rhoGal))\|_{L^\infty(0,T;\calH)}\,\le\,K_2^*\,.
\end{align}
\Ethm

\vspace{2mm}
\Brem
Notice that the pointwise condition \eqref{separ} is meaningful, since it follows from \cite[Sect.~8,~Cor.~4]{Simon}
and \eqref{regrho} that $\rho\in C^0(\overline Q)$ (and thus, in particular, that $\rhoG\in C^0(\overline\Sigma)$).
\Erem
\Brem
\pier{About \eqref{separ}, let us point out that,
u}nfortunately, we are unable to show a uniform in $\alpha\in (0,1]$ separation property. In fact, it may well
happen that, for $\alpha\searrow0$, we have $\,\rho_*(\alpha)\searrow -1\,$ and/or $\,\rho^*(\alpha)\nearrow +1$.
\Erem
\noindent
{\sc Proof of Theorem~2.3:} \quad The existence of a unique solution with the regularity \eqref{regmu} and \eqref{regrho},
which satisfies the separation property \eqref{separ}, is a direct consequence of \cite[Thm.~2.8]{CGS13}. 
In order to establish the global bound \eqref{wsb1}, we now follow the line of a priori estimates 
carried out in \cite{CGS13}, showing that the bounds derived there are in fact independent of $\alpha\in (0,1]$
in our special situation. In the following, we denote by $C$ positive constants that may depend on the data
of the system but neither on $u\in\uad$ nor on $\alpha\in (0,1]$. For the sake of a simpler \pier{notation}, 
we will also suppress the
superscript $\alpha$ in the calculations, writing it only at the end of each estimation. We also assume that
an arbitrary, but fixed, $u\in\uad$ is given. Observe that then \,$\|u\|_{\calX}\le R_0$, which will be
used repeatedly without further reference.

\vspace{2mm}\noindent
\underline{\sc First estimate:}

\vspace{1mm}\noindent Let $t\in (0,T]$ be arbitrary and $0<s\le t$. 
We insert $\,(v,\vG)=\Mu(s)\,$ in \eqref{was1} and $\,(v,\vG)=(\dt\rho,\dt\rhoG)(s)\,$ in \eqref{was2}, 
add the resulting equations, and integrate over $[0,t]$. Adding the expression \,$\,\intQt \rho\,\dt\rho
\,+\,\intSt \rhoG\,\dt\rhoG\,\,$ to both sides, we obtain the identity
\begin{align}
&\intQt|\nabla\mu|^2+\intSt|\nablaG\muG|^2 +\tauO\intQt|\dt\rho|^2+\tauG\intSt|\dt\rhoG|^2\non\\[1mm]
&+\,\frac 12\,\|\Rho(t)\|^2_{\calV}+\iO\vp(\alpha)h(\rho(t))+\iG\vp(\alpha)h(\rhoG(t))\non\\[1mm]
&=\,\frac 12\,\|(\rho_0,\rho_{0|\Gamma}\|_{\calV}^2+\iO\vp(\alpha)h(\rho_0)+\iG\vp(\alpha)h(\rho_{0|\Gamma})
+\intQt\rho \,u\cdot\nabla\mu\non\\[1mm]
&\quad\, +\intQt(\rho-{\pi}(\rho))\,\dt\rho+\intSt(\rhoG-{\piG}(\rhoG))\dt\rhoG\,,
\end{align}   
where, owing to the general assumptions, all of the terms on the \lhs\ are nonnegative and the first three terms on the \rhs\ are finite \pier{and uniformly bounded}. Now, recalling the separation property~\eqref{separ} \pier{and assumption~(A3)}, we conclude from Young's inequality that
the last two integrals on the \rhs\ are bounded by an expression of the form 
$\,\,C+\frac{\tauO}2\intQt|\dt\rho|^2+\frac{\tauG}2\intSt|\dt\rhoG|^2\,$. Moreover, owing to Young's inequality, \pier{we have that}
\begin{align}
&\intQt \rho\,u\cdot\nabla\mu\le\int_0^t\!\|\rho(s)\|_2\,\|u(s)\|_\infty\,\|\nabla \mu(s)\|_2\,ds
\,\le\,\frac 12\intQt|\nabla\mu|^2\,+\,C\,.
\end{align}
We thus  can infer from Gronwall's lemma the estimate
\begin{align}
\label{esti1}
&\|(\rhoal,\rhoGal)\|^2_{H^1(0,T;\calH)\cap L^\infty(0,T;\calV)} \,+\intQ|\nabla\mual|^2+\intS|\nablaG\muGal|^2\non\\[1mm]
&+\,\|\vp(\alpha)h(\rhoal)\|_{L^\infty(0,T;L^1(\Omega))}\,+\,\|\vp(\alpha)h(\rhoGal)\|_{L^\infty(0,T;L^1(\Gamma))}
\,\le\,C \quad\forall\,\alpha\in (0,1].
\end{align}

\vspace{2mm}\noindent
\underline{\sc Second estimate:}

\vspace{1mm}\noindent
Let $\hat m(t):={\rm mean}\,(\mu(t),\muG(t))$ for $t\in[0,T]$. Recalling \eqref{meanal}, we note that
$\,(v,\vG):=(\rho(t)-\hat r,\rhoG(t)-\hat r)\in {\cal V}_0\,$ for all $t\in [0,T]$. Inserting this in \eqref{was2},
where we temporarily omit the argument $\,t$, we obtain the identity
\begin{align}
\label{esti2}
&\iO\vp(\alpha)h'(\rho)(\rho-\hat r)\,+\iG\vp(\alpha)h'(\rhoG)(\rhoG-\hat r)\non\\[1mm]
\separa
&=\,-\tauO\iO\dt\rho(\rho-\hat r)\,-\,\tauG\iG\dt\rhoG(\rhoG-\hat r)\,-\iO|\nabla\rho|^2\,-
\iG|\nablaG\rhoG|^2\non\\[1mm]
&\quad\,\,-\iO {\pi}(\rho)(\rho-\hat r)\,-\iG {\piG}(\rhoG)(\rhoG-\hat r)\,+\iO(\mu-\hat m)(\rho-\hat r)\non\\[1mm]
&\quad\,\,+\iG(\muG-\hat m)(\rhoG-\hat r)\,.
\end{align}
\quad At this point, we recall that $-1<\hat r<1$. We thus may argue as in \cite[p.~908]{GMS} to conclude that there
exist constants $\,\delta_0>0\,$ and $\,C_0>0$, which do not depend on $\alpha\in (0,1]$, such that
$$
\vp(\alpha)h'(r)(r-\hat r)\,\ge\,\delta_0\,\vp(\alpha)|h'\pier{(r)}|-C_0 \quad\forall \,r\in (-1,1) \quad\forall\,
\alpha\in (0,1].
$$
\pier{Due to \eqref{separ},} the function  $\,\rho-\hat r\,$ is bounded on $\overline{Q}$\pier{;} we thus can infer from \eqref{esti2}, by just employing the Cauchy-Schwarz inequality, that
\begin{align}
\label{esti3}
&\delta_0\iO|\vp(\alpha)h'(\rho)|\,+\,\delta_0\iG |\vp(\alpha)h'(\rhoG)|\non\\[1mm]
&\le\,C\left(1\,+\,\juerg{\|\dt\rho\|_H+\|\dt\rhoG\|_{H_\Gamma}}
+\|\Rho\|_{\calV}^2+\juerg{\|(\mu-\hat m,\muG-\hat m)\|_{\calH}}\right)
\non\\[1mm]
&\le\,C\left(1\,+\,\juerg{\|\dt\rho\|_H+\|\dt\rhoG\|_{H_\Gamma}\,+\|(\mu-\hat m,\muG-\hat m)\|_{\calH}}\right),
\end{align}
where the last inequality follows from \eqref{esti1}. 
Now, we recall the definition \eqref{normaVz} and the fact that $\,\|\cdot\|_{{\cal V}_0}\,$ is equivalent
to the standard norm on $\calV_0$. Therefore,
$$
\|(\mu-\hat m,\muG-\hat m)\|_{\calH}\,\le\,C\,\|(\mu-\hat m,\muG-\hat m)\|_{{\calV}_0}
\,=\,C\,\|(\nabla\mu,\nablaG\muG)\|_{\calH}\,.
$$
Hence, combining this estimate with \eqref{esti1} and \eqref{esti3}, we can conclude that 
\begin{equation}
\label{esti4}
\|\vp(\alpha)h'(\rhoal)\|_{L^2(0,T;L^1(\Omega))}\,+\,\|\vp(\alpha)h'(\rhoGal)\|_{L^2(0,T;L^1(\Omega))}\,\le \,C
\quad\forall\,\alpha\in (0,1].
\end{equation}
At this point, we can insert $(v,\vG)=\pier{(1,1)}$ in \eqref{was2}, which then yields that the
function $\,t\mapsto {\rm mean}(\mual(t),\muGal(t))\,$ is bounded in $L^2(0,T)$, uniformly in
$\alpha\in (0,1]$. In view of \eqref{esti1}, we have thus shown that 
\begin{equation}
\label{esti5}
\|(\mual,\muGal)\|_{L^2(0,T;\calV)}\,\le\,C\quad\forall \,\alpha\in (0,1].
\end{equation}

\vspace{2mm}\noindent
\underline{\sc Third estimate:}

\vspace{1mm}\noindent
Next, we \pier{take} $(v,\vG)=(\vp(\alpha)h'(\rho(s)),\vp(\alpha)h'(\rhoG(s)))\in\calV$ in \eqref{was2},
where $0\le s\le t$ for some $t\in (0,T]$. Integrating over $[0,t]$, we obtain the identity
\begin{align}
&\tauO\iO\vp(\alpha)h(\rho(t))+\tauG\iG\vp(\alpha)h(\rhoG(t))+\intQt\vp(\alpha)h''(\rho)|\nabla\rho|^2
+\intSt\vp(\alpha)h''(\rhoG)|\nablaG\rhoG|^2\non\\[1mm]
&+\intQt|\vp(\alpha)h'(\rho)|^2+\intSt|\vp(\alpha)h'(\rhoG)|^2\non\\[1mm]&=\,\tauO\iO\vp(\alpha)h(\rhoz)+\tauG\iG\vp(\alpha)h(\rho_{0|\Gamma})+\intQt(\mu-{\pi}(\rho))\vp(\alpha)h'(\rho)\non\\[1mm]
&\quad\,\,+\intSt(\muG-{\piG}(\rhoG))\vp(\alpha)h'(\rhoG)\,,
\end{align} 
where (note that $h''\ge 0$) all of the terms on the \lhs\ are nonnegative and the first two summands on
the \rhs\ are bounded uniformly in $\alpha\in (0,1]$. Hence, in view of \eqref{esti5}, 
a simple application of
Young's inequality leads to the conclusion that
\begin{equation}
\label{esti6}
\|(\vp(\alpha)h'(\rhoal),\vp(\alpha)h'(\rhoGal))\|_{L^2(0,T;\calH)}\,\le\,C\quad\forall\,\alpha\in (0,1].
\end{equation} 
Direct comparison in \eqref{as2} then shows that also
\begin{equation}
\label{esti7}
\|\Delta\rhoal\|_{L^2(0,T;H)}\,\le\,C\quad\forall\,\alpha\in (0,1].
\end{equation}
Let us exploit \eqref{esti7}. Indeed, invoking \eqref{esti1}, \eqref{CO1} and \eqref{CO2},
we conclude that
\begin{equation}
\label{esti8}
\|\rhoal\|_{L^2(0,T;H^{3/2}(\Omega))}\,+\,\|\dn\rhoal\|_{L^2(0,T;\HG)}\,\le\,C
\quad\forall\, \alpha\in (0,1].
\end{equation}
Then comparison in \eqref{as4}, using \eqref{esti1}\pier{, \eqref{esti5}} and \eqref{esti6}, implies that
\begin{equation}
\label{esti9}
\|\DeltaG\rhoGal\|_{L^2(0,T;\HG)} \,\le\,C\quad\forall\,\alpha\in (0,1],
\end{equation}
and we conclude from \eqref{esti1}, \eqref{CO3} and \eqref{CO4} that
\begin{equation}
\label{esti10}
\|(\rhoal,\rhoGal)\|_{L^2(0,T;\calW)}\,\le\,C\quad\forall\,\alpha\in (0,1].
\end{equation}
Next, since $u\in\uad$, we readily infer from \eqref{as1} and \eqref{esti1}
that
\begin{equation}
\label{esti11}
\|\Delta\mual\|_{L^2(0,T;H)}\,\le\,C\quad\forall\,\alpha\in (0,1],
\end{equation}
whence, in view of \eqref{esti5}, \eqref{CO1} and \eqref{CO2},
\begin{equation}
\label{esti12}
\|\mual\|_{L^2(0,T;H^{3/2}(\Omega))}\,+\,\|\dn\mual\|_{L^2(0,T;\HG)}
\,\le\,C\quad\forall\alpha\in (0,1].
\end{equation}
Hence, by virtue of \eqref{as3} and \eqref{esti1}, \pier{we have that}
\begin{equation}
\label{esti13}
\|\DeltaG\muGal\|_{L^2(0,T;\HG)}\,\le\,C \quad\forall\,\alpha\in (0,1],
\end{equation}
and we can argue as above to arrive at the estimate
\begin{equation}
\label{esti14}
\|(\mual,\muGal)\|_{L^2(0,T;\calW)} \,\le\,C\quad\forall\,\alpha\in (0,1].
\end{equation}

\vspace{2mm}\noindent
\underline{\sc Fourth estimate:}
\par\nobreak
\vspace{1mm}\noindent
We now argue formally, noting that the following arguments can be made rigorous by,
e.g., using finite differences in time. At first, we note that $\,\mean \dt\Rho=0$ a.e.
in $(0,T)$, by \pier{\eqref{meanal}}. 
Hence, $\Xi{(t)}:=\calN(\dt\Rho{(t)})\in\calV_0$ is well defined {for a.a.~$t\in(0,T)$}. 
Now, we differentiate both \eqref{was1}
and \eqref{was2} (formally) with respect to time, test the resulting identities by $\Xi$ and $\dt\Rho$, respectively, 
and add the results. Now observe that, by \eqref{perdefN} and
\eqref{defN},
\begin{align*}
\intQt\nabla\dt\mu\cdot\nabla\xi\,+\intSt\nablaG\dt\muG\cdot\nablaG\xiG
\,=\,\intQt \dt\mu\,\dt\rho\,+\intSt\dt\muG\,\dt\rhoG\,.
\end{align*}  
Hence,
recalling \eqref{propNdta}, and integrating the expressions containing $u$ (formally) by parts, 
 we  arrive at the identity
\begin{align}\label{esti15} 
&\frac 12\,\|\dt\Rho(t)\|_*^2\,+\,\frac{\tauO}2\iO|\dt\rho(t)|^2\,+\,\frac{\tauG}2\iG|\dt\rhoG(t)|^2
+\intQt|\nabla\dt\rho|^2+\intSt|\nablaG\dt\rhoG|^2\non\\[1mm]
&+\intQt\vp(\alpha)h''(\rho)|\dt\rho|^2+\intSt\vp(\alpha)h''(\rhoG)|\dt\rhoG|^2\non\\[1mm]
&=\,I_0+\intQt\nabla\dt\rho\cdot u\xi+\intQt\nabla \rho\cdot\dt u\,\xi-\intQt {\pi'}(\rho)|\dt\rho|^2
-\intSt {\piG'}(\rhoG)|\dt\rhoG|^2\,,
\end{align}
where
\begin{equation}\label{esti16}
I_0\,:=\,\frac 12\,\|\dt\Rho(0)\|_*^2\,+\,\frac{\tauO}2\iO|\dt\rho(0)|^2\,+\,\frac{\tauG}2\iG|\dt\rhoG(0)|^2\,.
\end{equation}
Noting that $\vp(\alpha)h''\ge 0$, we may omit the two nonnegative summands in the second line of \eqref{esti15}, and
thus obtain an inequality which has exactly the same form as the inequality \cite[Eq.~(7.1)]{CGS13}. 
We thus may repeat 
the estimates carried out in \cite{CGS13} in order to conclude that (cf., \cite[Eq.~(7.3)]{CGS13})
\begin{equation}\label{esti17}
\|(\rhoal,\rhoGal)\|_{W^{1,\infty}(0,T;\calH)\cap H^1(0,T;\calV)}\,\le\,C\quad\forall\,\alpha\in(0,1]
\end{equation}
\pier{and actually realize that the estimate is uniform with respect to $\alpha$.}
 
\vspace{2mm}\noindent
\underline{\sc Fifth estimate:}

\vspace{1mm}\noindent
Recalling that $\hat m(t):={\rm mean}\,(\mu(t),\muG(t))$ for $t\in[0,T]$, we test  \eqref{was1} by the 
$\,\calV_0$--valued function $\,\Mu-\hat m(1,1)$. Using the fact that the norm \eqref{normaVz} is equivalent
to the standard norm on $\calV_0$, we obtain, for almost every $t\in (0,T)$,
\begin{align}\label{esti18}
&\iO|\nabla\mu|^2+\iG|\nablaG\muG|^2\,=\,-\iO\dt\rho(\mu-\hat m)-\iG\dt\rhoG(\muG-\hat m)-\iO\rho u\cdot\nabla\mu
\non\\[2mm]
&\le \,C\,\|\dt\Rho\|_{\calH}\,\|\Mu-\hat m(1,1)\|_{{\cal V}_0}\,+\,\|u\|_{\pier{L^\infty(Q)}}\,\|\rho\|_2\,\|\nabla\mu\|_2
\non\\[2mm]
&\le \,C\left(\|\nabla\mu\|_2+\|\nablaG\muG\|_2\right)\,.
\end{align}
Consequently, \pier{we deduce that}
\begin{align}
\label{esti19}
&\|\nabla\mual\|_{L^\infty(0,T;H)}\,+\,\|\nablaG\pier{\muGal}\|_{L^\infty(0,T;H)}\,\le\,C\quad\forall\,\alpha\in (0,1],\\[2mm]
\label{esti20}
&\|{(\mual-\hat m,\muGal-\hat m)}\|_{L^\infty(0,T;\cal V)}\,\le\,C\quad\forall\,\alpha\in (0,1]. 
\end{align}

\vspace{2mm}\noindent
\underline{\sc Sixth estimate:}

\vspace{1mm}\noindent
At first, we directly obtain from \eqref{esti3}, \eqref{esti17}, and \eqref{esti20}, that
\begin{equation}
\label{esti21}
\|\vp(\alpha)h'(\rhoal)\|_{L^\infty(0,T;L^1(\Omega))}+\|\vp(\alpha)h'(\rhoGal)\|_{L^\infty(0,T;L^1(\Omega))}
\,\le\,C\quad\forall\alpha\in (0,1].
\end{equation}
Therefore, if we take $\,(v,\vG)=(1,1)/(|\Omega|+|\Gamma|)$\, in \eqref{was2},
\pier{for almost every $t\in (0,T)$ we infer that}
\begin{align}
|{\rm mean}\,\Mu(t)|\,&\le\,C\,\|\dt\Rho\|_{L^\infty(0,T;\calH)}\,+\,C\,\|\vp(\alpha)h'(\rho)+{\pi}(\rho)\|_
{L^\infty(0,T;L^1(\Omega))}\non\\[1mm]
&\quad\,+\,C\,\|\vp(\alpha)h'(\rhoG)+{\piG}(\rhoG)\|_{L^\infty(0,T;L^1(\Omega))}\,\le\,C\,.
\end{align} 
By virtue of \eqref{esti20}, this shows that
\begin{align}
\label{esti22}
\|(\mual,\muGal)\|_{L^\infty(0,T;\calV)}\,\le\,C\quad\forall\,\alpha\in (0,1].
\end{align} 

At this point, we observe that \eqref{as1}, \eqref{esti1}, \eqref{esti17}, and the fact that $u\in\uad$,
imply that
\begin{equation}
\|\Delta\mual\|_{L^\infty(0,T;H)}\,\le\,\|\dt\rhoal\|_{L^\infty(0,T;H)}+ \|u\cdot\nabla\rhoal\|_{L^\infty(0,T;H)}\,\le\,C
\quad\forall\,\alpha\in (0,1].
\end{equation}
In view of \eqref{esti22}, we are therefore in the same situation as in the third estimation above after the
proof of \eqref{esti11} (only that we have $L^\infty$ with respect to time in place of $L^2$). We thus may 
argue as in the estimates \eqref{esti12}--\eqref{esti14} to conclude that
\begin{equation}
\label{esti23}
\|(\mual,\muGal)\|_{L^\infty(0,T;\calW)}\,\le\,C\quad\forall\,\alpha\in (0,1].  
\end{equation}

\vspace{2mm}\noindent
\underline{\sc Seventh estimate:}

\vspace{1mm}\noindent
Finally, we insert $\,(v,\vG)=(\vp(\alpha)h'(\rho),\vp(\alpha)h'(\rhoG))\,$ in \eqref{was2}. Employing the
estimates shown previously, we readily obtain that, almost everywhere on $(0,T)$,
\begin{align}
\label{esti24}
&\iO|\vp(\alpha)h'(\rho)|^2+\iG|\vp(\alpha)h'(\rhoG)|^2+\iO\vp(\alpha)h''(\rho)|\nabla\rho|^2
+\iG\vp(\alpha)h''(\rhoG)|\nablaG\rhoG|^2\non\\[1mm]
&=\,\iO\vp(\alpha)h'(\rho)\left(-\tauO\dt\rho+\mu-{\pi}(\rho)\right)+\iG\vp(\alpha)h'(\rhoG)\left(
-\tauG \dt\rhoG+\muG-{\piG}(\rhoG)\right)\non\\[1mm]
&\le\,C\,+\,\frac 12\iO|\vp(\alpha)h'(\rho)|^2\,+\,\frac 12\iG|\vp(\alpha)h'(\rhoG)|^2\,.
\end{align}
Consequently, \pier{we have that}
\begin{align}\label{esti25}
\|(\vp(\alpha)h'(\rhoal),\vp(\alpha)h'(\rhoGal))\|_{L^\infty(0,T;\calH)}\,\le\,C\quad\forall\,\alpha\in (0,1].
\end{align}
But then, by virtue of \eqref{as2} and the previous estimates, \pier{it is clear that}
\begin{equation*}
\|\Delta\rhoal\|_{\pier{L^\infty(0,T;H)}}\,\le\,C\quad\forall\,\alpha\in (0,1],
\end{equation*}
and, arguing as in the third estimate, we infer that
\begin{equation}\label{esti26}
\|(\rhoal,\rhoGal)\|_{L^\infty(0,T;\calW)}\,\le\,C\quad\forall\alpha\in (0,1],
\end{equation}
which concludes the proof of the assertion.\QED

\vspace{5mm}\noindent
\Brem
By virtue of the well-posedness result given by Theorem~2.3, 
the control-to-state operator $\calS_\alpha: u\mapsto (\Mual,\Rhoal)$ is well defined as a 
mapping between $\uad\subset\calX$ and the space defined by the regularity stated in \eqref{regmu}, \eqref{regrho}. 
In particular, this also holds true for its second component $\,\calS^2_\alpha:u\mapsto\Rhoal$.
\Erem


\section{Existence and approximation of optimal controls}
\setcounter{equation}{0}

In this section, we aim to approximate optimal pairs of (${\cal P}_0$). To this end, we consider for $\alpha\in (0,1]$ the optimal control problem

\vspace{5mm} \noindent
(${\cal P}_\alpha$) \quad Minimize the cost functional
$\,{\cal J}((\rhoal,\rhoGal),u)\,$ for $\,u\in\uad$,  subject 
to the state system   
\hspace*{12mm} \eqref{as1}--\eqref{as5}.

\vspace{5mm} \noindent Assuming generally that (A1)--(A5) are fulfilled, we obtain from
\cite[Thm.~4.1]{CGS14} that this optimal control problem has an
optimal pair $(((\mual,\muGal),(\rhoal,\rhoGal)),\ual)$, for every $\alpha\in (0,1]$.
Our first aim in this section is to prove the following approximation result:

\Bthm
Suppose that the assumptions {\rm (A1)}--{\rm (A5)}, \eqref{defh} and \eqref{phiat0} 
are satisfied, and let sequences $\{\alpha_n\}\subset (0,1]$ and 
$\{\un\}\subset \uad$ be given such that $\alpha_n\searrow 0$ and
$\un\to u$ weakly-star in ${\cal X}$ for some $u\in\uad$. Then there is a subsequence
$\{\alpha_{n_k}\}$ of $\{\alpha_n\}$ such that for $k\to\infty$ it holds,
with $\,(\Mun,\Rhon):={\cal S}_{\alpha_n}(\un)$, $n\in\nz$,  
\begin{align}
\label{conmu}
&(\mu^{\alpha_{n_k}},\muG^{\alpha_{n_k}})\to \Mu\quad\mbox{\em weakly-star in }\,
 L^\infty(0,T;\calW) , \\[1mm]
\label{conrho}
&(\rho^{\alpha_{n_k}},\rhoG^{\alpha_{n_k}})\to \Rho\quad\mbox{\em weakly-star in }\,W^{1,\infty}(0,T;\calH)\cap H^1(0,T;\calV)
\cap L^\infty(0,T;\calW),\\[1mm]
\label{conxi}
&(\vp(\alpha_{n_k})\,h'(\rho^{\alpha_{n_k}}),\vp(\alpha_{n_k})\,h'(\rhoG^{\alpha_{n_k}}))\to\Xi\quad\mbox{\em weakly-star in }\,L^\infty(0,T;\calH),
\end{align}
where $(\Mu, \Rho,\Xi)$ is a  solution to the state
system {\rm \eqref{ss1}--\eqref{ss7}} associated with $u$. Moreover, \eqref{conrho} holds true for the
entire sequence $\{\alpha_n\}$.  Finally, with
$\,{\cal S}_0^2(u):=\Rho\,$ it holds that
\begin{align}
\label{conj1}
{\cal J}({\cal S}_0^2(u),u)\,&\le\,\liminf_{n\to\infty}\,{\cal J}({\cal S}_{\alpha_n}^2
(\un),\un),\\[1mm]
\label{conj2}
{\cal J}({\cal S}_0^2(v),v)\,&=\lim_{n\to\infty}\,{\cal J}({\cal S}_
{\alpha_n}^2(v),v)\quad\forall\,v\in\uad.
\end{align}
\Ethm
 
\vspace{3mm}\noindent
{\sc Proof:} \,\quad
Let $\,\{\alpha_n\}\subset (0,1]\,$ be any sequence such that $\alpha_n\searrow 0$ as $n\to\infty$, and suppose that $\un\to u$ weakly-star in ${\cal X}$ for some $u\in\uad$. By virtue of Theorem~2.3, there are a subsequence of $\{\alpha_n\}$, which is again indexed by $n$, and three pairs $\,\Mu,\Rho,\Xi\,$ such that
the convergence results \eqref{conmu}--\eqref{conxi} hold true. Moreover, from standard compact 
embedding results (cf. \cite[Sect.~8, Cor.~4]{Simon}) we can infer that
\begin{align}
\label{strcon1}
\rhon&\to\rho\quad\mbox{strongly in }\,L^2(0,T;V)\cap C^0(\overline Q),
\end{align}
which also yields that
\begin{equation}
\label{strcon2}
\rhoGn\to\rhoG\quad\mbox{strongly in }\,C^0(\overline \Sigma)\,.
\end{equation}
In particular, $(\rho(0),\rhoG(0))=(\rho_0,\rho_{0|\Gamma})$ and $\,\rhoG=\rho_{|\Sigma}$. In addition, 
we obviously have that
\begin{align}
\label{strcon3}
{\pi}(\rhon)&\to {\pi}(\rho)\quad\mbox{strongly in $\,C^0(\overline Q)$},\\[1mm]
\label{strcon4}
{\piG}(\rhoGn)&\to{\piG}(\rhoG)\quad\mbox{strongly in $\,C^0(\overline \Sigma)$}.
\end{align} 
Moreover, it is easily verified that, at least weakly in $L^1(Q)$,
\begin{align}
&\nabla\rhon\cdot\un\to \nabla\rho\cdot u\,.
\end{align}

Combining the above convergence results, we may pass to the limit
as $n\to\infty$ in the equations (\ref{as1})--(\ref{as5}) (written for 
$\alpha=\alpha_n$ and $u=\un$) to find that $(\Mu,\Rho,\linebreak \Xi)$ and $\,u\,$
satisfy the equations \eqref{ss1}, (\ref{ss2}), \eqref{ss4},
\eqref{ss5}, and \eqref{ss7}.  Thus, in order to show that 
 $(\Mu,\Rho,\Xi)$ is in fact a solution to the
problem (\ref{ss1})--(\ref{ss7}) corresponding to $u$, 
it remains to
 show that $\,\xi\in\partial I_{[-1,1]}(\rho)\,$ a.\,e. in $Q$
and $\,\xiG\in \partial I_{[-1,1]}(\rhoG)\, $ a.\,e. in $\Sigma$. 
To this end, recall that $\,h\,$ is convex and bounded in $[-1,1]$ and
\juerg{that} both $\,h\,$ and $\,\vp\,$ are nonnegative.
We thus have, for every $n\in\enne$,
\begin{align}
\label{t3.1}
&0\,\leq\, 
\intQ\varphi(\alpha_n)\,h(\rho^{\alpha_n})\,\le\,\intQ\varphi(\alpha_n)\,h(z)\,+
\,\intQ\varphi(\alpha_n)\,h'( \rho^{\alpha_n})\, (\rho^{\alpha_n} - z)
\nonumber\\[1mm]
&\mbox{for all }\,z\in {\cal K}:=\{v\in {L^2(Q)}:|v|\leq 1\text{ a.e. in }Q\}\,. 
\end{align}
Thanks to (\ref{phiat0}), the first two integrals tend to zero as $n\to\infty$. Hence, invoking  
(\ref{conxi}) and \juerg{(\ref{strcon1})},  the passage to the limit as $n\to\infty$ yields
\begin{equation}
\label{t3.2}
\intQ\xi\,( \rho -z)\geq 0\quad\forall z\in {\mathcal{K}} . 
\end{equation}
Inequality (\ref{t3.2}) entails that $\xi$ is an element of the subdifferential of the extension $\calI $ of $ I_{[-1,1]}$ to $L^2(Q)$, which means that $\xi \in \partial \,\calI (\rho)$ or, equivalently (cf.~\cite[Ex.~2.3.3., p.~25]{Brezis}),  that
$\xi\in\partial I_{[-1,1]}(\rho)$ {a.\,e. in $Q$}. Similarly, we can prove that 
$\xi_\Gamma\in\partial I_{[-1,1]}(\rho_\Gamma)$ a.\,e. in $\Sigma$.

\vspace{1mm}
\quad We have thus shown that, for a suitable subsequence of $\{\alpha_n\}$,
we have the convergence properties \eqref{conmu}--\eqref{conxi}, where $(\Mu,\Rho,\Xi)$ is a solution to the state system \eqref{ss1}--\eqref{ss7}. But, according to Theorem~2.2, the component $\,\Rho\,$ is the same for any such solution. 
This entails that the  convergence
properties \eqref{conrho}, \eqref{strcon1}--\eqref{strcon4} are in fact valid for the entire sequence
$\{\alpha_n\}$. This finishes the proof of the
first claims of the theorem.

It remains to show the validity of \eqref{conj1} and \eqref{conj2}. 
In view of \eqref{conrho}, the inequality \eqref{conj1}  is an immediate consequence of the 
weak and weak-star sequential semicontinuity properties of the cost functional ${\cal J}$. To establish
the identity \eqref{conj2}, let $v\in \uad$ be arbitrary  and  put
$(\rho^{\alpha_n},\rhoG^{\alpha_n})={\cal S}_{\alpha_n}^2(v)$, for $n\in\nz$. 
Taking Theorem~2.3 into account, and arguing as in the first part of this proof,
we can conclude that $\{\SALN(v)\}$ converges to $\Rho={\calS}_0^2(v)$ in the sense
of \eqref{conrho}. In particular, we have (recall \eqref{strcon1} and \eqref{strcon2})
$$\SALN(v)\to {\cal S}_0^2(v)\quad\mbox{strongly in }\,C^0(\overline Q)\times C^0(\overline\Sigma).
$$
As the cost functional ${\cal J}$ is obviously continuous in the variables $(\rho,\rhoG)$
with respect to the strong topology of $\,C^0(\overline Q)\times C^0(\overline\Sigma)$,
we may thus infer that \eqref{conj2} is valid. \QED

\vspace{4mm}
\Bcor
Under the assumptions of Theorem~3.1, the optimal control problem {\rm (${\cal P}_0$)}  has a least one solution.
\Ecor

\noindent
{\sc Proof:} \quad Pick an arbitrary sequence $\{\alpha_n\}$ such that $\alpha_n\searrow0$ as $n\to\infty$.  
Then, by virtue of \cite[Thm.~4.1]{CGS14}, the optimal control problem (${\cal P}_{\alpha_n}$) has for
every $n\in\nz$ an optimal pair $((\Rhon,\Mun),\un)$, where $\,\un\in\uad\,$ and  $\,\Rhon=
\SALN(\un)$. Since $\uad$ is a bounded subset of ${\cal X}$, we may without loss of generality assume
that $\,\un\to u\,$ weakly-star in ${\cal X}$ for some $\,u\in\uad$. 
Then, for some solution $(\Mu,\Rho,\Xi)$ to the state system \eqref{ss1}--\eqref{ss7} associated
with $u$, we conclude from Theorem~3.1 the convergence properties \eqref{conrho}, \eqref{strcon1},
\eqref{strcon2}, and \eqref{conj2}.
Invoking the optimality of $((\Mun,\Rhon),\un)$ for (${\cal P}_{\alpha_n}$),
we then find, for every $v\in\uad$, that 
\begin{align}
\label{tr3.3}
&{\cal J}(\Rho,u)\,=\,{\cal J}({\cal S}_{0}^2(u),u)\,\le\,
\liminf_{n\to\infty}\,{\cal J}({\mathcal S}_{\alpha_n}^2(u^{\alpha_n}),u^{\alpha_n})
 \nonumber\\[1mm]
&\leq\,\liminf_{n\to\infty}\,{\cal J}({\cal S}_{\alpha_n}^2(v),v)\, =\,\lim_{n\to\infty} {\cal J}({\cal S}_{\alpha_n}^2(v),v)\,=\,
{\cal J}({\mathcal S}_{0}^2(v),v),
\end{align}  
which yields that $u$ is an optimal control for (${\cal P}_0$) with the associate state
$(\Mu,\Rho, \linebreak \Xi)$. The assertion is thus proved.
\QED

\vspace{7mm}\quad
Corollary~3.2 does not yield any information on whether every solution to the optimal control problem $({\mathcal{P}}_{0})$ can be approximated by a sequence of solutions to the problems $({\mathcal{P}}_{\alpha})$. 
As already announced in the Introduction, we are not able to prove such a 
general `global' result. Instead, we 
can only give a `local' answer for every individual optimizer of $({\mathcal{P}}_{0})$. For this purpose,
we employ a trick due to Barbu~\cite{Barbu}. To this end, let $\bar u\in\uad$
be an arbitrary optimal control for $({\mathcal{P}}_{0})$, and let $((\bar\mu,\bar\mu_\Gamma),
(\bar\rho,\bar\rho_\Gamma), (\bar\xi,\bar\xi_\Gamma))$ be any associated solution  to the state system (\ref{ss1})--(\ref{ss7}) in the sense of 
Theorem~\pier{2.2}. In particular, $\,(\bar\rho,\bar\rho_\Gamma)={\cal S}_0^2 (\bar u)$. We associate with this 
optimal control the {\em adapted cost functional}
\begin{equation}
\label{cost2}
\widetilde{\cal J}((\rho,\rhoG),u):={\cal J}((\rho,\rhoG),u)\,+\,\frac{1}{2}\,\|u-\bar{u}\|^2_{(L^2(Q))^3}
\end{equation}
and a corresponding {\em adapted optimal control problem},

\vspace{4mm}\noindent
($\widetilde{\mathcal{P}}_{\alpha}$)\quad Minimize $\,\, \widetilde {\cal J}(\Rho,u)\,\,$
for $\,u\in\uad$, subject to the condition that  
(\ref{as1})--(\ref{as5}) \hspace*{11mm} be satisfied.

\vspace{3mm}
With a standard direct argument that needs no repetition here, we can show the following 
result.

\Blem
Suppose that the assumptions {\rm (A1)}--{\rm (A5)}, {\rm (\ref{defh})} and {\rm (\ref{phiat0})} are
satisfied, and let $\alpha\in (0,1]$. Then the optimal control problem 
$(\widetilde{\cal P}_\alpha)$ admits a solution.
\Elem
  
\vspace{3mm}
We are now in the position to give a partial answer to the question raised above. We have the following result.

\Bthm
Let the assumptions {\sc (A1)}--{\sc (A5)}, {\rm (\ref{defh})} and {\rm (\ref{phiat0})} be fulfilled, suppose that 
$\bar u\in \uad$ is an arbitrary optimal control of {\rm $({\mathcal{P}}_{0})$} with associated state  
$((\bar\mu,\bar\mu_\Gamma),(\bar\rho,\bar\rho_\Gamma),\linebreak(\bar\xi,\bar\xi_\Gamma))$,
and let $\,\{\alpha_n\}\subset (0,1]\,$ be
any sequence such that $\,\alpha_n\searrow 0\,$ as $\,n\to\infty$. Then there exist a 
subsequence $\{\alpha_{n_k}\}$ of $\{\alpha_n\}$, and, for every $k\in\nz$, an optimal control
 $\,u^{\alpha_{n_k}}\in \uad\,$ of the adapted problem {\rm $(\widetilde{\mathcal{P}}_{\alpha_{n_k}})$}
 with associated state $((\mu^{\alpha_{n_k}},\mu_\Gamma^{\alpha_{n_k}}),(\rho^{\alpha_{n_k}},
 \rho_\Gamma^{\alpha_{n_k}}))$ such that, as $k\to\infty$,
\begin{align}
\label{tr3.4}
&u^{\alpha_{n_k}}\to \bar u\quad\mbox{strongly in }\,(L^2(Q))^3,
\end{align}
and such that the properties {\rm \eqref{conmu}--\eqref{conxi}} are satisfied, where $(\mu,\muG),(\rho,
\rhoG),(\xi,\xiG)$ are replaced by  $(\bar\mu,\bar\mu_\Gamma),(\bar\rho,\bar\rho_\Gamma),
(\bar\xi,\bar\xi_\Gamma)$. Moreover, we have 
\begin{align}
\label{tr3.5}
&\lim_{k\to\infty}\,\widetilde{{\cal J}}((\rho^{\alpha_{n_k}},\rho_\Gamma^{\alpha_{n_k}}),
u^{\alpha_{n_k}})\,=\,  {\cal J}((\bar\rho,\bar\rho_\Gamma),\bar u)\,.
\end{align}
\Ethm

\vspace{1mm}\noindent
{\sc Proof:} \quad\,Let $\alpha_n \searrow 0$ as $n\to\infty$. For any $ n\in\nz$, we pick an optimal control 
$u^{\alpha_n} \in \uad\,$ for the adapted problem $(\widetilde{\cal P}_{\alpha_n})$ and denote by 
$((\mu^{\alpha_n},\muGn),(\rhon,\rhoGn))$ the associated solution to the problem (\ref{as1})--(\ref{as5}) for $\alpha=\alpha_n$ and $u=\un$. 
By the boundedness of $\uad$ in $\calX$, there is some subsequence $\{\alpha_{n_k}\}$ of $\{\alpha_n\}$ such that
\begin{equation}
\label{ugam}
u^{\alpha_{n_k}}\to u\quad\mbox{weakly-star in }\,{\cal X}
\quad\mbox{as }\,k\to\infty,
\end{equation}
with some $u\in\uad$. Thanks to Theorem~3.1, the convergence properties \eqref{conmu}--\eqref{conxi}
hold true, where $(\Mu,\Rho,\Xi)$ is some solution to the state system
\eqref{ss1}--\eqref{ss7}. In particular,  $((\Mu,\Rho,\Xi),u)$
is admissible for (${\cal P}_0$). 

\vspace{2mm}\quad
We now aim to prove that $u=\bar u$. Once this is shown, then the uniqueness result of 
Theorem~2.2 yields that also $(\rho,\rhoG)=(\bar\rho,\bar\rho_\Gamma)$, 
which implies that 
the properties {\rm \eqref{conmu}--\eqref{conxi}} are satisfied, where $(\mu,\muG),(\rho,
\rhoG),(\xi,\xiG)$ are replaced by  $(\bar\mu,\bar\mu_\Gamma),(\bar\rho,\bar\rho_\Gamma),
(\bar\xi,\bar\xi_\Gamma)$.

Now observe that, owing to the weak sequential lower semicontinuity of 
$\widetilde {\cal J}$, 
and in view of the optimality property of $((\bar\mu,\bar\mu_\Gamma),(\bar\rho,\bar\rho_\Gamma),
(\bar\xi,\bar\xi_\Gamma),\bar u)$ for problem $({\cal P}_0)$,
\begin{align}
\label{tr3.6}
&\liminf_{k\to\infty}\, \widetilde{\cal J}((\rho^{\alpha_{n_k}},
\rho_\Gamma^{\alpha_{n_k}}),u_\Gamma^{\alpha_{n_k}})
\ge \,{\cal J}((\rho,\rhoG),u)\,+\,\frac{1}{2}\,
\|u-\bar{u}\|^2_{(L^2(Q))^3}\nonumber\\[1mm]
&\geq \, {\cal J}((\bar\rho,\bar\rho_\Gamma),\bar u)\,+\,\frac{1}{2}\,\|u-\bar{u}\|^2_{(L^2(Q))^3}\,.
\end{align}
On the other hand, the optimality property of  $\,(((\mu^{\alpha_{n_k}},\mu_\Gamma^{\alpha_{n_k}}),
(\rho^{\alpha_{n_k}},\rho_\Gamma^{\alpha_{n_k}})),u^{\alpha_{n_k}})
\,$ for problem $(\widetilde {\cal P}_{\alpha_{n_k}})$ yields that
for any $k\in\nz$ we have
\begin{equation}
\label{tr3.7}
\widetilde {\cal J}((\rho^{\alpha_{n_k}},\rho_\Gamma^{\alpha_{n_k}}),u^{\alpha_{n_k}})\, =\,
\widetilde {\cal J}({\cal S}^2_{\alpha_{n_k}}(u^{\alpha_{n_k}}),
u^{\alpha_{n_k}})\,\le\,\widetilde {\cal J}({\cal S}^2_{\alpha_{n_k}}
(\bar u),\bar u)\,,
\end{equation}
whence, taking the limit superior as $k\to\infty$ on both sides and invoking (\ref{conj2}) in
Theorem~3.1,
\begin{align}
\label{tr3.8}
&\limsup_{k\to\infty}\,\widetilde {\cal J}((\rho^{\alpha_{n_k}},
\rho_\Gamma^{\alpha_{n_k}}),u^{\alpha_{n_k}})
\nonumber\\[1mm]
&\le\,\widetilde {\cal J}({\calS}_0^2(\bar u),\bar u) 
\,=\,\widetilde {\cal J}((\bar\rho,\bar\rho_\Gamma),\bar u)
\,=\,{\cal J}((\bar\rho,\bar\rho_\Gamma),\bar u)\,.
\end{align}
Combining (\ref{tr3.6}) with (\ref{tr3.8}), we have thus shown that 
$\,\frac{1}{2}\,\|u-\bar{u}\|^2_{(L^2(Q))^3}=0$\,,
so that $\,u=\bar u\,$  and thus also $\,(\rho,\rhoG)
=(\bar\rho,\bar\rho_\Gamma)$. 
Moreover, (\ref{tr3.6}) and (\ref{tr3.8}) also imply that
\begin{align}
\label{tr3.9}
&{\cal J}((\bar\rho,\bar\rho_\Gamma),\bar u) \, =\,\widetilde{\cal J}((\bar\rho,\bar\rho_\Gamma),\bar u)
\,=\,\liminf_{k\to\infty}\, \widetilde{\cal J}((\rho^{\alpha_{n_k}},
\rho_\Gamma^{\alpha_{n_k}}), u^{\alpha_{n_k}})\nonumber\\[1mm]
&\,=\,\limsup_{k\to\infty}\,\widetilde{\cal J}((\rho^{\alpha_{n_k}},
\rho_\Gamma^{\alpha_{n_k}}), u^{\alpha_{n_k}}) \,
=\,\lim_{k\to\infty}\, \widetilde{\cal J}((\rho^{\alpha_{n_k}},
\rho_\Gamma^{\alpha_{n_k}}), u^{\alpha_{n_k}})\,,
\end{align}                                     
which proves {(\ref{tr3.5})} and, at the same time, also (\ref{tr3.4}). This concludes the proof
of the assertion.\QED

\section{The optimality system}
\setcounter{equation}{0}
In this section, we aim to establish first-order necessary optimality conditions for the optimal control problem $({\mathcal{P}}_{0})$.  This will be achieved by a passage to the limit as $\alpha\searrow 0$ in the first-order necessary optimality conditions for the adapted optimal control problems $(\widetilde{\mathcal{P}}_{\alpha})$ that can by derived as in \cite{CGS14} with only
minor and obvious changes. This procedure will yield certain generalized first-order
necessary optimality conditions in the limit. In this entire section, we generally assume that
$\,h\,$ is given by (\ref{defh}) and that (\ref{phiat0}) and
the assumptions {\sc (A1)}--{\sc (A5)} are
satisfied. \pier{In addition, we} assume that the following \elvis{condition is} fulfilled: 

\vspace{2mm}\noindent
\elvis{(A6) \quad  $\tauO=\tauG=:\tau>0$.}

\vspace{2mm}\noindent
We also assume that a fixed optimal control $\bar u\in \uad$ for 
$({\cal P}_0)$
is given, along
with a corresponding solution $((\bar\mu,\bar\mu_\Gamma),(\bar\rho,\bar\rho_\Gamma),(\bar\xi,\bar\xi_\Gamma))$ to the 
state system (\ref{ss1})--(\ref{ss7}) in the sense of Theorem~2.2. That is, we have 
$(\bar\rho,\bar\rho_\Gamma)={\calS}_0^2(\bar u)$, as well as $\bar\xi\in\partial I_{[-1,1]}(\bar\rho)$ a.\,e. in $Q$
and $\bar\xi_\Gamma\in \partial I_{[-1,1]}(\bar\rho_\Gamma)$ a.\,e. on $\Sigma$.

\vspace{2mm}
 We begin our analysis by formulating the adjoint state system for the adapted
control problem $(\widetilde{\mathcal{P}}_{\alpha})$ corresponding to $\bar u$.
To this end, let us assume that, for some $\alpha\in (0,1]$, $u^\alpha\in\uad$ is an arbitrary optimal control for 
$(\widetilde{\mathcal{P}}_{\alpha})$ and that $((\mual,\muGal),(\rhoal,\rhoGal))$ 
is the (unique) solution to the associated state system (\ref{as1})--(\ref{as5}). In particular, 
$((\mual,\muGal),(\rhoal,\rhoGal))={\cal S}_\alpha(u^\alpha)$,  the solution enjoys the 
regularity properties (\ref{regmu}) and (\ref{regrho}), and it satisfies the global bounds
\eqref{wsb1} and the separation property \eqref{separ}. The associated adjoint system has the following 
variational form (cf., \cite[Eqs.~(4.7)--(4.9)]{CGS14}):
\begin{align}
&-\left\langle\dt\left(\pal+\tau\qal,\pGal+\tau\qGal\right),(v,\vG)\right\rangle_{\calV}\,+
\iO\nabla\qal\cdot\nabla v\,+\iG\nablaG\qGal\cdot\nablaG\vG
\non \\[1mm]
&\quad{}+\iO(\vp(\alpha)h''(\rhoal)+ {\pi'}(\rhoal))\,\qal\,v\,+
\iG(\vp(\alpha)h''(\rhoGal)+{\piG'}(\rhoGal))\,\qGal\,\vG\,-\iO\ual\cdot\nabla\pal\,v
\non\\[1mm]
&=\,\iO\beta_1(\rhoal-\widehat\rho_Q)\,v \,+\iG\beta_2(\rhoGal-\widehat\rho_\Sigma)\,\vG
\quad \mbox{\,a.e. in $(0,T)$}, \quad\forall\, (v,\vG)\in\calV,
\label{adj1}
\\[2mm]
\separa
&\iO\nabla\pal\cdot\nabla v\,+\iG\nablaG\pGal\cdot\nablaG\vG\,=\,\iO\qal v\,+\iG\qGal \vG 
\non\\[1mm]
&\quad \mbox{\,a.e. in $(0,T)$}, \quad\forall\, (v,\vG)\in\calV,
\label{adj2}
\\[2mm]
&\left\langle\left(\pal+\tau\qal,\pGal+\tau\qGal\right)(T),(v,\vG)\right\rangle_{\calV}
\,=\,\iO\beta_3(\rhoal(T)-\widehat\rho_\Omega)v\,+\iG\elvis{\beta_4}
(\rhoGal(T)-\widehat\rho_\Gamma)\vG
\non\\[1mm]
&\quad\forall\,(v,\vG)\in {\calV},
\label{adj3}
\end{align}

\vspace{2mm}\noindent
which corresponds to the {backward} problem
\begin{align}
&-\,\dt\left(\pal+\tau\qal\right)-\Delta\qal+\vp(\alpha)h''(\rhoal)\qal+{\pi'}(\rhoal)\qal-\ual\cdot\nabla\pal=\beta_1(\rhoal-\widehat\rho_Q)\non\\[1mm]
&\quad\,\mbox{ and }\,-\Delta\pal=\qal\,\mbox{ in }\,Q,\\[2mm]
&-\,\dt\left(\pGal+\tau\qGal\right)+\dn\qal-\DeltaG\qGal+\vp(\alpha)h''(\rhoGal)\qGal+{\piG'}(\rhoGal)\qGal=\beta_2(\rhoGal-\widehat\rho_\Sigma),\non\\[1mm]
&\quad\, \dn\pGal-\DeltaG\pGal=\qGal,\quad \, p^\alpha_{|\Sigma}=\pGal\,\,\mbox{ and }\,\, 
q^\alpha_{|\Sigma}=\qGal\,\quad\mbox{on \,$\,\Sigma$},\\[2mm]
\label{ini}
&\quad\left(\pal+\tau\qal,\pGal+\tau\qGal\right)(T)=\left(\beta_3(\rhoal(T)-\widehat\rho_\Omega),
\elvis{\beta_4}(\rhoGal(T)-\widehat\rho_\Gamma)\right)\,.
\end{align}

\vspace{2mm}\noindent
According to \cite[Thm.~4.4]{CGS14}, the adjoint system \eqref{adj1}--\eqref{adj3} enjoys for every $\alpha\in (0,1]$
a unique solution $\left((\pal,\pGal),(\qal,\qGal)\right)$ such that
\begin{align}
\label{regadj1}
&\left(\pal,\pGal\right)\in \pier{L^\infty(0,T;\calV)}, \quad \left(\qal,\qGal\right)\in L^\infty(0,T;\calH)\cap L^2(0,T;\calV),\\[1mm]
\label{regadj2}
&\left(\pal+\tau\qal,\pGal+\tau\qGal\right)\in H^1(0,T;{\calV}^*).
\end{align}

\vspace{1mm}\noindent
Observe that, owing to  \eqref{regadj1} and \eqref{regadj2}, 
$$(\pal+\tau\qal,\pGal+\tau\qGal)\in (H^1(0,T;\calV^*)\cap L^2(0,T;\calV))\subset
C^0([0,T];\calH),
$$
by continuous embedding. In particular, the {final} condition \eqref{adj3} is in fact satisfied 
in the form \eqref{ini}. 
Moreover, arguing as in the derivation of \cite[Thm.~4.6]{CGS14}, we can infer that for 
any such solution $\left((\pal,\pGal),(\qal,\qGal)\right)$ there holds the variational 
inequality
\begin{align}
\label{vugal}
\intQ \left(\rho^\alpha\nabla\pal+\beta_5\ual+(\ual-\bar u)\right)\cdot  (v-\ual)\,\ge\,0
\quad\forall\,v\in\uad\,.
\end{align}

We now try to find bounds that are uniform with respect to
$\alpha\in (0,1]$. 
To this end, we define for $\alpha\in (0,1]$ the quantities 
\begin{align}
\label{quant1}
&\vp^\alpha_Q:=\beta_1(\rhoal-\widehat\rho_Q),\mbox{ } 
\vp^\alpha_\Sigma:=\beta_2(\rhoGal-\widehat\rho_\Sigma),\mbox{ }
\vp^\alpha_\Omega:=\beta_3(\rhoal(T)-\widehat\rho_\Omega),\mbox{ }
\vp^\alpha_\Gamma:=\elvis{\beta_4}(\rhoGal(T)-\widehat\rho_\Gamma),
\end{align}
noting that (A3), (A4) and \eqref{wsb1} imply that
\begin{align}\label{quant2}
&\|\vp^\alpha_Q\|_{L^2(Q)}\,+\,\|\vp^\alpha_\Sigma\|_{L^2(\Sigma)}
\,+\,\|\vp^\alpha_\Omega\|_{L^2(\Omega)}\,+\,\|\vp^\alpha_\Gamma\|_{L^2(\Gamma)}\non\\[1mm]
&+\,\|{\pi'}(\rhoal)\|_{L^\infty(Q)}\,+\,\|{\piG'}(\rhoGal)\|_{L^\infty(\Sigma)}
\,\le\,K_3^* \quad\forall\,\alpha\in (0,1],
\end{align}
with a constant $K_3^*>0$ that depends only on the data of the system. 

In view of the low regularity of the adjoint state variables, the derivation of 
uniform bounds makes it necessary to argue 
by approximation, following an idea introduced in the proof of \cite[Thm.~4.4]{CGS14}.  
Namely, for fixed $\alpha\in (0,1]$, we approximate $(\vp_\Omega^\alpha,\vp_\Gamma^\alpha)$ by pairs 
$(\vp_\Omega^{\alpha,\eps},\vp_\Gamma^{\alpha,\eps})$, $\eps\in (0,1]$, which satisfy
\begin{align}
\label{t41.1}
&(\vp_\Omega^{\alpha,\eps}/\tau,\vp_\Gamma^{\alpha,\eps}/\tau)\in \calV,\quad
(\vp_\Omega^{\alpha,\eps},\vp_\Gamma^{\alpha,\eps})\to (\vp^\alpha_\Omega,
\vp^\alpha_\Gamma)\quad\mbox{in $\,\calH$ \,as \,$\eps\to 0$,}
\end{align}
and consider for every $\eps\in (0,1]$ the approximating system
\begin{align}
\label{app1}
&-\iO\dt(\palep+\tau \qalep)v\,-\iG\dt(\pGalep+\tau \qGalep)\vG
\,+\iO\nabla\qalep\cdot\nabla v\,+\iG\nablaG\qGalep\cdot\nablaG\vG\non\\
&+\iO(\vp(\alpha)h''(\rhoal)+{\pi'}(\rhoal))\qalep\,v \,+
\iG(\vp(\alpha)h''(\rhoGal)+{\piG'}(\rhoGal))\qGalep\,\vG\non\\
&-\iO\ual\cdot\nabla\palep\,v \,=\,\iO\vp_Q^\alpha \,v+\iG\juerg{\vp^\alpha_\Sigma}\,\vG
\qquad\mbox{$\forall\,(v,\vG)\in\calV$ and a.e. in $(0,T)$},\\[2mm]
\label{app2}
&-\eps\iO\dt\palep\,v\,-\eps\iG\dt\pGalep\,\vG\,+\iO\nabla\palep\cdot\nabla v
\,+\iG\nablaG\pGalep\cdot\nablaG\vG\non\\
&=\,\iO\qalep\,v\,+\iG\qGalep\,\vG
\qquad\mbox{$\forall\,(v,\vG)\in\calV$ and a.e. in $(0,T)$},\\[2mm]
\label{app3}
&(\palep,\pGalep)(T)\,=\,(0,0),\quad\,(\qalep,\qGalep)(T)\,=\,
(\vp_\Omega^{\alpha,\eps}/\tau,\vp_\Gamma^{\alpha,\eps}/\tau)\,.
\end{align}

According to \cite[Thm.~4.3]{CGS14}, the system \eqref{app1}--\eqref{app3}
enjoys for every $\eps\in (0,1]$ a unique solution $\,((\palep,\pGalep),(\qalep,\qGalep))\,$ such that
\begin{equation}\label{t41.2}
(\palep,\pGalep), \,(\qalep,\qGalep)\in H^1(0,T;\calH)\cap L^\infty(0,T;\calV).
\end{equation}
Moreover, it was shown in the proof of \cite[Thm.~4.4]{CGS14} that there is some sequence $\eps_n\searrow0$ such that,
as $n\to\infty$,
\begin{align}
\label{eps1}
&(p^{\alpha,\eps_n},p_\Gamma^{\alpha,\eps_n})\to (\pal,\pGal)\quad\mbox{weakly\pier{-star} in }\,\pier{L^\infty(0,T;\calV)},\\[1mm]
\label{eps2}
&(q^{\alpha,\eps_n},q_\Gamma^{\alpha,\eps_n})\to (\qal,\qGal)\quad\mbox{weakly-star in }\,
L^\infty(0,T;\calH)\cap L^2(0,T;\calV),\\[1mm]
\label{eps3}
&\dt(p^{\alpha,\eps_n}+\tau q^{\alpha,\eps_n},p_\Gamma^{\alpha,\eps_n}+\tau q_\Gamma^{\alpha,\eps_n})\to \dt(\pal+\tau\qal,\pGal+\tau\qGal)\quad\mbox{weakly in }\,L^2(0,T;\calV^*),\\[1mm]
\label{eps4}
&\eps_n\,\dt(p^{\alpha,\eps_n},p_\Gamma^{\alpha,\eps_n})\to (0,0)\quad\mbox{strongly in }\,L^2(0,T;\calH),
\end{align}
where $((\pal,\pGal),(\qal,\qGal))$ is the solution to the adjoint system  \eqref{adj1}--\eqref{adj3} having the
regularity properties \eqref{regadj1}--\eqref{regadj2}. Notice that \eqref{eps1}--\eqref{eps3} imply that also
\begin{align}\label{eps5}
(p^{\alpha,\eps_n}+\tau q^{\alpha,\eps_n},p_\Gamma^{\alpha,\eps_n}+\tau q_\Gamma^{\alpha,\eps_n})
\to (\pal+\tau\qal,\pGal+\tau\qGal) \non\qquad\\
\pier{\mbox{strongly in }\,C^0([0,T];\calV^*)
\cap L^2(0,T;\calH)},
\end{align}
so that the Cauchy condition \eqref{adj3} is meaningful. In the following, we will always work with the
particular sequence $\{\eps_n\}$.

Next, we establish
uniform bounds for the approximating solutions. To simplify {the} notation,
we omit in the following estimate the superscript $^{\alpha,\eps}$, writing it only
at the end of the respective estimations. We also recall the definition 
of $Q^t$ and $\Sigma^t$, for $t\in [0,T)$, given in \eqref{defQt}, and we denote
by $C_i$, $i\in\enne$, positive constants that may depend on the data, but neither on
$\eps\in (0,1]$ nor on $\alpha\in (0,1]$.
 
We test \eqref{app1} by $(q,\qG)$, integrate over~$(t,T)$, and
account for the Cauchy conditions \eqref{app3}, to obtain the identity
\Bsist
  && - \bintQt \dt p \, q
  - \bintSt \dt\pG \, \qG
  + \frac \tau 2 \iO |q(t)|^2
  + \frac \tau 2 \iG |\qG(t)|^2
  + \bintQt |\nabla q|^2
  + \bintSt |\nablaG\qG|^2
  \non
  \\
  &&+ \bintQt \varphi(\alpha)h''(\rhoal)|q|^2 +\bintSt \varphi(\alpha)h''(\rhoGal)|\qGal|^2\non\\
  && = \frac \tau 2 \iO |\phi_\Omega^{\alpha,\eps}/\tau|^2
  + \frac \tau 2 \iG |\phi_\Gamma^{\alpha,\eps}/\tau|^2 \,+ \bintQt \ual \cdot \nabla p \, q
  - \bintQt {\pi'}(\rhoal)|q|^2
  - \bintSt {\piG'}(\rhoGal) |\qG|^2
  \non
  \\
  \label{equa1}
  && \quad {}
    + \bintQt \phi_Q^\alpha\, q
  + \bintSt \phi_\Sigma^\alpha \, \qG \,.
\Esist
At the same time, we test \eqref{app2} by $-\dt(p,\pG)$ and integrate over~$(t,T)$ 
to obtain the identity
\Beq
  \eps \bintQt |\dt p|^2
  + \eps \bintSt |\dt\pG|^2
  + \frac 12 \iO |\nabla p(t)|^2
  + \frac 12 \iG |\nablaG\pG(t)|^2
  = - \bintQt q \dt p
  - \bintSt \qG \dt\pG \,.\label{equa2}
 \Eeq
Now, we add \eqref{equa1} and \eqref{equa2}, observing that four terms cancel out
and that the two summands in the second line of \eqref{equa1} are nonnegative. Omitting these two
summands and the first two summands in the \lhs\ of \eqref{equa2}, we then arrive at the
inequality
\begin{align}
&\frac \tau 2 \iO |q(t)|^2
  + \frac \tau 2 \iG |\qG(t)|^2
  + \bintQt |\nabla q|^2
  + \bintSt |\nablaG\qG|^2
+\frac 12\iO|\nabla p(t)|^2+\frac 12\iG|\nablaG \pG(t)|^2\non\\
&\le\,C_1\,+\,C_2\Big(\bintQt|q|^2+\bintSt|\qG|^2\Big)\,+\,\bintQt \ual \cdot \nabla p \, q\,,
\end{align}
where we have used \eqref{quant2}, \eqref{t41.1} and Young's inequality. Now, by Young's inequality, and
since $\ual\in\uad$, 
\Bsist
  && \bintQt \ual \cdot \nabla p \, q
  \,\leq\, \norma{\pier{u^\alpha}}_{\pier{L^\infty(Q)}} \int_t^T \|\nabla p(s)\|_2 \, \|q(s)\|_2 \, ds
  \non
  \\
  && \hspace*{25mm}\le\,\bintQt |\nabla p|^2\,+\,C_3\bintQt | q|^2\,. 
\Esist
Therefore, invoking Gronwall's lemma, we can infer that
\begin{align}\label{estaleps}
&\|(q^{\alpha,\eps},q_\Gamma^{\alpha,\eps})\|_{L^\infty(0,T;\calH)\cap L^2(0,T;\calV)}
\,+\,\supess_{t\in (0,T)} \Big(\iO|\nabla p^{\alpha,\eps}(t)|^2 +\iG|\nablaG p_\Gamma^{\alpha,\eps}(t)|^2\Big)
\,\le\,C_4
\end{align}
for all $\alpha\in (0,1]$ and $\eps\in (0,1]$. We thus can conclude from the weak and weak-star sequential lower
semicontinuity of norms, taking the limit as $\eps_n\searrow 0$, that 
\Beq
\|(\qal,\qGal)\|_{\L\infty\calH\cap\L2\calV} \,+\,
   \supess_{t\in (0,T)}\Big(\iO|\nabla\pal(t)|^2 \,+\,\iG|\nablaG\pGal(t)|^2\Big)\,\le\,C_4
  \label{estal}
\Eeq
for all $\alpha\in (0,1]$.

\Brem\label{remark}
In the proof of \cite[Thm.~4.4]{CGS14}, further estimates for the approximations $((p^{\alpha,\eps},p_\Gamma
^{\alpha,\eps}),(q^{\alpha,\eps},q_\Gamma^{\alpha,\eps}))$ could be derived. However, a closer look at these estimations
reveals that the resulting bounds depend on the special choice of $\alpha\in (0,1]$ and may become infinite as
$\alpha\searrow 0$. 
In particular, 
while it is clear that 
\Beq\label{meanq}
{\rm mean}\,(\qal(t),\qGal(t))=0\quad\mbox{for all \,$t\in [0,T]$\, and \,$\alpha\in (0,1]$,}   
\Eeq
as one immediately sees by inserting $(v,\vG)=(1,1)$ in \eqref{adj2},
it seems to be impossible to derive a uniform bound for the mean value
of $(\pal,\pGal)$, the main reason being that the separation constants $\rho_*(\alpha),\rho^*(\alpha)$ 
introduced in Theorem~2.3, which were implicitly used in the argument to control the
expressions $\varphi(\alpha)h''(\rhoal)\qal$ and $\varphi(\alpha)h''(\rhoGal)\qGal$, may approach $\pm 1$ as $\alpha\searrow0$. The difficulty becomes apparent if we observe that insertion of $(v,\vG)=(1,1)$ in \eqref{adj1} and integration of the resulting identity over
$[t,T]$, where $t\in [0,T]$, yields the representation formula 
{(by~also owing to~\eqref{meanq}) 
\begin{align}
\label{meanp}
&{\rm mean}\,(\pal(t),\pGal(t))
=\,\frac 1{|\Omega|+|\Gamma|}\Bigg[
   {-}\bintQt(\varphi(\alpha)h''(\rhoal(t))+{\pi'}(\rhoal(t)))\qal(t)\non\\
& {}-\bintSt(\varphi(\alpha)h''(\rhoGal(t))+{\piG'}(\rhoGal(t)))\qGal(t)
\non\\
&+\iO\beta_3(\rhoal(T)-\widehat\rho_\Omega) + \iG\elvis{\beta_4}(\rhoGal(T)-\widehat\rho_\Gamma)
+\bintQt\beta_1(\rhoal-\widehat\rho_Q)+\bintSt\beta_2(\rhoGal-\widehat\rho_\Sigma)\Bigg]\,.
\end{align}
}%
\Erem

In order to be able to derive a meaningful adjoint system for problem (${\cal P}_0$), we 
thus have to eliminate the
mean value of $(\pal,\pGal)$ from the problem, thereby avoiding the difficulty mentioned above. To this
end, we follow a strategy introduced in \cite{CGS3} and \cite{CFGS1}: 
{by recalling~\eqref{meanq}, it}
 follows from \eqref{adj2} and the definition \eqref{defN} of the operator ${\cal N}$ the identity
\Bsist
&& (\pal(t),\pGal(t))-{\mean}\,(\pal(t),\pGal(t)){(1,1)}\,=\,{\cal N}(\qal(t),\qGal(t))
\non
\\
&& {\quad\mbox{for all \,$t\in[0,T]$\, and \,$\alpha\in (0,1]$}\,.}
\label{p=Nq}
\Esist
Since $(\qal,\qGal)$ is uniformly bounded in $L^\infty(0,T;\calH)$, in particular, we can infer from
\cite[Lem.~3.1]{CGS13} that $(\xi^\alpha,\xi_\Gamma^\alpha):={\cal N}(\qal,\qGal)$ belongs to 
$ L^\infty(0,T;\calW\cap\calHz)$,
solves the boundary value problem
$$
-\Delta \xi^\alpha(t)=\qal(t) \quad\mbox{a.e. in $\Omega$}, \qquad \dn \xi^\alpha(t)-\DeltaG 
\xi_\Gamma^\alpha(t)=\qGal(t)\quad\mbox{a.e. on $\Gamma$},
$$  
for almost every $t\in (0,T)$, and satisfies the uniform bound
\Beq\label{regNq}
\|{\cal N}(\qal,\qGal)\|_{L^\infty(0,T;\calW)}\,\le\,C_5\quad\forall\,\alpha\in (0,1].  
\Eeq

Now recall that $\calV=\calVz\oplus\Span\{(1,1)\}$, where $\calVz$ is defined in \eqref{defcalVz}. Notice also that,
by virtue of \pier{\cite[Lem.~5.1~and~Cor.~5.3]{CGS3}}, it holds that $\VG=\{\vG:(v,\vG)\in \calVz\}$ and that 
${\cal H}_0$ is dense in $\calVz$. We thus can construct the Hilbert triple $\,\calVz\subset\calHz
\subset \calVz^*\,$ with dense and compact embeddings, that is, we identify $\calHz$ with a subspace
of $\calVz^*$ in such a way that 
\begin{equation}\label{unafestasuiprati}
\langle(w,w_\Gamma),(v,\vG)\rangle_{\calVz}=\iO w\,v+\iG w_\Gamma\,\vG\quad\forall\,(w,w_\Gamma)
\in \calHz,\quad\forall\,(v,\vG)\in \calVz\,.
\end{equation} 
Notice that the embedding $\,(H^1(0,T;\calVz^*)\cap L^2(0,T;\calVz))\subset
C^0([0,T];\calHz)$ is continuous. Observe also that, because of the zero mean value condition,
the first components \,$v\,$ of the elements $(v,\vG)\in\calVz$ do not span the whole space
$C_0^\infty(\Omega)$, so that variational equalities with test functions from $\calVz$
cannot directly be interpreted as equations in the sense of distributions.

At this point, the additional \pier{assumption (A6) comes} into play. To this end,
recall that $\pier{(z^\alpha, z^\alpha_\Gamma)}:=\dt(\pal+\tau\qal,\pGal+\tau\qGal)$ belongs to $L^2(0,T;\calV^*)$ and
thus also to $L^2(0,T;\calVz^*)$. We now aim to show a global bound for the
family $\,\{\pier{(z^\alpha, z^\alpha_\Gamma)}\}_{\alpha\in(0,1]}\,$ that will prove to be fundamental for the
subsequent argumentation. To this end, we introduce the spaces
\begin{align}
\label{defZ}
&{\cal Z}:=(H^1(0,T;V^*)\times H^1(0,T;V_\Gamma^*))\cap L^2(0,T;\calVz)\,,
\\[2mm]
\label{defZ0}
&{\cal Z}_0:=\{(v,\vG)\in{\cal Z}:\,(v(0),\vG(0))=(0,0)\}\,,
\end{align}
which are Banach spaces when endowed with the natural norm of ${\cal Z}$. Moreover,
${\cal Z}$ is continuously embedded in $\,C^0([0,T];\calHz)$, so that the initial
condition encoded in \eqref{defZ0} is meaningful. In addition,  ${\cal Z}_0$ is a closed subspace
of $Y\times Y_\Gamma$, where 
\begin{align}
\label{defY}
&Y:=H^1(0,T; V^*)\cap L^2(0,T;V)
\aand
Y_\Gamma:=H^1(0,T; V_\Gamma^*)\cap L^2(0,T;\VG)
\end{align}
are Banach spaces when endowed with their natural norms. 
\pier{It} then follows (cf., e.g., \cite[Prop.~2.6]{CGS1}) that 
the elements ${F\in{\mathcal Z}^*_0}$ are exactly those that are of the form 
\begin{align}
\label{eq:4.15}
\left\langle F,(\eta,\eta_\Gamma)\right\rangle_{{\cal Z}_0} \,=\,
\left\langle z, \eta  \right\rangle_Y\,+\,
\left\langle z_\Gamma,\eta_\Gamma \right\rangle_{Y_\Gamma} \quad 
\hbox{ for all $(\eta,\eta_\Gamma)\in {\mathcal Z}_0$,}
\end{align}
with some $z\in Y^*$ and $z_\Gamma\in Y_\Gamma^*$. 
{Thus, we can write
\Bsist
  \< F , (\eta,\eta_\Gamma) >_{\calZ_0}
  = \ioT \< z(t) , \eta(t) >_V \, dt
  + \ioT \< z_\Gamma(t) , \eta_\Gamma(t) >_{\VG} \, dt 
  \quad \hbox{for every $(\eta,\eta_\Gamma)\in\calZ_0$}.
  \non
\Esist 
Moreover, even though the pair $(z,z_\Gamma)$ associated with $F\in{\mathcal Z}^*_0$ is not unique,
the above representation formula allows} 
us to give a proper meaning to statements like 
$$
(z^\alpha , z_\Gamma^\alpha ) \to  (z , z_\Gamma )
\quad \hbox{ weakly in } {\mathcal Z}^*_0 .  
$$

Now let $(v,\vG)\in {\cal Z}_0$ be arbitrary. Then $(v,\vG)(0)=(0,0)$, ${\rm mean\,}(v(t),\vG(t))=0$ for all
$t\in [0,T]$, and ${\rm mean\,}(\dt(v,\vG)(t))=0$ for almost every $t\in (0,T)$. 
{Thus, from one side, we have by~\eqref{p=Nq} \juerg{that}
\Bsist
  && \< \dt(v,\vG) , (p^\alpha,p_\Gamma^\alpha) >_{\calV}
  = \< \dt(v,\vG) , (p^\alpha,p_\Gamma^\alpha) - \mean(p^\alpha,p_\Gamma^\alpha)(1,1) >_{\calV}
  \non
  \\
  && = \< \dt(v,\vG) , \calN(q^\alpha,q_\Gamma^\alpha) >_{\calV}
  \quad \aet .
  \label{gia}
\Esist
On the other hand,} 
using \eqref{adj3}
and the fact \juerg{that} both $(v,\vG)$ and $(\pal+\tau\qal,\pGal+\tau\qGal)$ belong to $H^1(0,T;\calV^*)\cap
L^2(0,T;\calV)$, 
{we see that}
\begin{align}\label{intpart}
&\int_0^T\langle -\pier{(z^\alpha, z^\alpha_\Gamma)}(t),(v,\vG)(t)\rangle_{\calV}\,dt\,=\,-\iO\beta_3(\rhoal(T)-\widehat\rho_\Omega)v(T)
-\iG\elvis{\beta_4(\rhoGal}(T)-\widehat\rho_\Gamma)\vG(T)\non\\[1mm] 
&+\int_0^T\langle\dt(v,\vG)(t),(\pal+\tau\qal,\pGal+\tau\qGal)(t)\rangle_{\calV}\,dt \,=:\,A^\alpha\,.
\end{align}
\pier{Thanks to \eqref{quant1}--\eqref{quant2}, the}
sum of the first two summands on the \rhs\ of \eqref{intpart}, which we denote by~$A_1^\alpha$, \elvis{satisfies the estimate}
\begin{align}\label{yesteryou}
|A_1^\alpha|\,&\le\,\|(\pier{\varphi^\alpha_\Omega},\pier{\varphi^\alpha_\Gamma})\|_{\elvis{\calH}}\,
\|(v,\vG)(T)\|_{\elvis{\calH}}\non\\[1mm]
&\le\,C_6\,\|(v,\vG)\|_{{\cal Z}_0} \quad\forall\,\alpha\in (0,1]\,,
\end{align}
where \pier{the continuity of the embedding $(H^1(0,T;\calVz^*)\cap L^2(0,T;\calVz))
\subset C^0([0,T];\calHz)$ has been used}. 
Moreover, the third summand on the \rhs\ of \eqref{intpart}, which we denote
by~$A_2^\alpha$, satisfies, in view of~{\eqref{gia}}, the identity
\begin{equation}\label{yesterday}
A_2^\alpha\,=\,\int_0^T\langle \dt(v,\vG)(t),({\cal N}(\qal(t),\qGal(t))+\tau(\qal(t),\qGal(t))\rangle_{\calVz}\,dt\,.
\end{equation}
In addition, \pier{from \eqref{estal} and \eqref{regNq} it follows} that
\begin{align}\label{azzurro}
|A_2^\alpha|\,&\le\,\int_0^T\|\dt(v,\vG)(t)\|_{\calVz^*}\,\|{\cal N}(\qal(t),\qGal(t))+\tau(\qal(t),\qGal(t))\|_{\calVz}\,dt
\non\\[1mm]
&\le\,C_7\,\|\dt(v,\vG)\|_{L^2(0,T;\calVz^*)}\,\le\,C_7\,\|(v,\vG)\|_{{\cal Z}_0}\quad\forall\,\alpha\in (0,1].
\end{align}
{}From the above estimates, we can infer that
\begin{align}\label{dtpal}
\|\dt(\pal+\tau\qal,\pGal+\tau\qGal)\|_{{\cal Z}_0^*}\,\le\,C_8\quad\forall\,\alpha\in(0,1],
\end{align}
and comparison in \eqref{adj1} yields that also
\begin{equation}\label{lamal}
\|(\varphi(\alpha)h''(\rhoal)\qal,\varphi(\alpha)h''(\rhoGal)\qGal)\|_{{\cal Z}_0^*}\,\le\,C_9\quad
\forall\,\alpha\in(0,1].
\end{equation}

We  are now in a position to state the first-order necessary optimality conditions for problem (${\cal P}_0$).

\Bthm
Suppose that {\rm (A1)--\elvis{(A6)}}, \eqref{defh}, \eqref{phiat0} are satisfied, and let $\bar u\in\uad$ be an optimal
control with associated state $((\bar\mu,\bar\mu_\Gamma),(\bar\rho,\bar\rho_\Gamma),(\bar\xi,\bar\xi_\Gamma))$ in the
sense of Theorem~2.2. Then there exist $(q,\qG),\,\eta,\,(\lambda,\lambda_\Gamma)$ such that the following
statements hold true: 

\vspace{1mm}\noindent
{\rm (i)} \hspace*{5mm} $(q,\qG)\in L^\infty(0,T;\calHz)\cap L^2(0,T;\calVz)$, \,\,${\cal N}(q,\qG)\in L^\infty(0,T;\calW
\cap \calHz)$, \,\,$(\lambda,\lambda_\Gamma)\in {\cal Z}_0^*$,\\[1mm]
\hspace*{11.5mm} and \,$\eta\in (L^\infty(0,T;H))^3$.

\vspace{1mm}
\noindent
{\rm (ii)} \hspace*{3.8mm} Adjoint system:
\begin{align}\label{adjoint}
&\langle(\lambda,\lambda_\Gamma),(v,\vG)\rangle_{{\calZ}_0}\,+\int_0^T\langle\dt(v,\vG)(t),
{\cal N}(q(t),\qG(t))+\tau(q(t),\qG(t))\rangle_{\calVz}\,dt\non\\[1mm]
&+\intQ {\pi'}(\bar\rho)q\, v\,+ \intS {\piG'}(\bar\rho_\Gamma)\qG\,\vG \,+\intQ \bar u\cdot\eta\,v \non\\[1mm]
&=\,\intQ\beta_1(\bar\rho-\widehat\rho_Q)v\,+\intS\beta_2(\bar\rho_\Gamma-\widehat\rho_\Sigma)\vG\,+
\iO\beta_3(\bar\rho(T)-\widehat\rho_\Omega)v(T)\non\\[1mm]
&\quad + \,\iG\elvis{\beta_4}(\bar\rho_\Gamma(T)-\bar\rho_\Gamma)\vG(T) \qquad\,\forall\,(v,\vG)\in {\cal Z}_0\,.
\end{align}

\vspace{1mm}\noindent
{\rm (iii)} \hspace*{2.9mm} Necessary optimality condition:
\begin{equation}\label{vug}
\intQ (\bar\rho\,\eta+\beta_5\,\bar u)\cdot(v-\bar u)\,\ge\,0\quad\forall\,v\in\uad\,.
\end{equation}

\Ethm 
	
\vspace{3mm}\noindent
{\sc Proof:} \quad  We pick a sequence $\{\alpha_n\}_{n\in\enne}\subset (0,1]$ such that $\alpha_n\searrow0$
and (cf. Theorem~3.1, Theorem~3.4, and \eqref{strcon1}--\eqref{strcon4})
\begin{align}
\label{conny1}
&u^{\alpha_n}\to \bar u \quad\mbox{strongly in \,$(L^2(Q))^3$},\\[2mm]
&(\rho^{\alpha_n},\rho_\Gamma^{\alpha_n})\to (\bar\rho,\bar\rho_\Gamma) \quad
\mbox{weakly-star in }\,W^{1,\infty}(0,T;\calH)\cap H^1(0,T;\calV)\cap L^\infty(0,T;\calW)\non\\
&\hspace*{39mm}\mbox{and strongly in }\,L^2(0,T;\calV)\cap (C^0(\overline Q)\times C^0(\overline\Sigma)),
\label{conny2}\\[2mm]
\label{conny3}
&({\pi'}(\rho^{\alpha_n}),{\piG'}(\rho_\Gamma^{\alpha_n}))\to ({\pi'}(\bar\rho),{\piG'}(\bar\rho_\Gamma))\quad
\mbox{strongly in }\,C^0(\overline Q)\times C^0(\overline\Sigma)\,.
\end{align}  
Moreover, in view of the estimates \eqref{estal}, \eqref{regNq}, and \eqref{lamal}, we may assume
that there are $(q,\qG)\in L^\infty(0,T;\calHz)\cap L^2(0,T;\calVz)$, $\eta\in (L^\infty(0,T;H))^3$,
and $(\lambda,\lambda_\Gamma)\in {\cal Z}_0^*$, such that
\begin{align}\label{conny4}
&(q^{\alpha_n},\qG^{\alpha_n})\to (q,\qG)\quad\,\mbox{weakly-star in }\,L^\infty(0,T;\calHz)\cap
L^2(0,T;\calVz),\\[2mm]
\label{conny5}
&\nabla p^{\alpha_n}\to\eta\quad\,\mbox{weakly-star in }\,(L^\infty(0,T;H))^3,\\[2mm]
\label{conny6}
&(\varphi(\alpha_n)h''(\rho^{\alpha_n})q^{\alpha_n},\varphi(\alpha_n)h''(\rho_\Gamma^{\alpha_n})q_\Gamma^{\alpha_n})
\to (\lambda,\lambda_\Gamma)\quad\,\mbox{weakly in }\, {\cal Z}_0^*,\\[2mm]
\label{conny7}
&{\cal N}(q^{\alpha_n},\qG^{\alpha_n})\to {\cal N}(q,\qG)\quad\,\mbox{weakly-star in }\,
L^\infty(0,T;\calW\cap\calHz).
\end{align}
Now we can take advantage of the identities \pier{\eqref{intpart} and \eqref{yesterday}}.
Indeed, if we restrict ourselves to test functions $(v,\vG)\in {\cal Z}_0$ and invoke the convergence
results \eqref{conny2}--\eqref{conny7}, then we may pass to the limit as $n\to\infty$ in the
equations \eqref{adj1}--\eqref{adj3} (written for $\alpha=\alpha_n$) to arrive at the conclusion that
we have the identity
\begin{align}
\label{yeah?}
&\langle(\lambda,\lambda_\Gamma),(v,\vG)\rangle_{{\calZ}_0}\,+\int_0^T\langle\dt(v,\vG)(t),
{\cal N}(q(t),\qG(t))+\tau(q(t),\qG(t))\rangle_{\calVz}\,dt\non\\[1mm]
&+\intQ {\pi'}(\bar\rho)q\, v\,+ \intS {\piG'}(\bar\rho_\Gamma)\qG\,\vG \,+\,
\lim_{n\to\infty}\intQ u^{\alpha_n}\cdot\nabla p^{\alpha_n}\,v \non\\[1mm]
&=\,\intQ\beta_1(\bar\rho-\widehat\rho_Q)v\,+\intS\beta_2(\bar\rho_\Gamma-\widehat\rho_\Sigma)\vG\,+
\iO\beta_3(\bar\rho(T)-\widehat\rho_\Omega)v(T)\non\\[1mm]
&\quad + \,\iG\elvis{\beta_4}(\bar\rho_\Gamma(T)-\bar\rho_\Gamma)\vG(T) \qquad\,\forall\,(v,\vG)\in {\cal Z}_0\,.
\end{align}
Therefore, in order to prove the validity of \eqref{adjoint}, we need to show that 
\begin{equation}\label{yeah!!}
\lim_{n\to\infty}\intQ u^{\alpha_n}\cdot\nabla p^{\alpha_n}\,v \,=\,\intQ\bar u\cdot\eta v
\quad\forall\,(v,\vG)\in {\cal Z}_0.
\end{equation} 
To this end, it suffices to establish the result for all test functions from   the set \,$\widetilde 
{\cal Z}_0:=\{(v,\vG)\in {\cal Z}_0:
v\in {L^2(0,T;C^0(\overline\Omega))}\}$. Indeed, since $\widetilde{\cal Z}_0$  is a dense subset of ${\cal Z}_0$,
\eqref{yeah!!} then follows from a simple density argument. Now let $(v,\vG)\in\widetilde{\cal Z}_0$. We 
have that 
\begin{align}
\intQ(u^{\alpha_n}\cdot\nabla p^{\alpha_n}-\bar u\cdot\eta)\,v \,=\,\intQ(u^{\alpha_n}-\bar u)\cdot\nabla
p^{\alpha_n}\,v\,+\intQ\bar u\cdot(\nabla p^{\alpha_n}-\eta)\,v\,.
\end{align}
Since $\bar u\in\uad$, we can infer from \eqref{conny5} that the second integral on the \rhs\ approaches zero as
$n\to\infty$. Moreover, we obtain from \eqref{conny1}, using \eqref{estal} and
H\"older's inequality, that
\begin{align}\label{yeah!!!}
&\Big|\intQ(u^{\alpha_n}-\bar u)\cdot\nabla p^{\alpha_n}\,v\Big|\,\le\,\int_0^T\|u^{\alpha_n}(t)-\bar u(t)\|_2
\,\|\nabla p^{\alpha_n}(t)\|_2\,\|v(t)\|_\infty\,dt\non\\[1mm]
&\le\,{\,\|u^{\alpha_n}-\bar u\|_{(L^2(Q))^3}\,\norma{\nabla p^{\alpha_n}}_{(\L\infty\Ldue)^3}\,\|v\|_{L^2(0,T;C^0(\overline\Omega))}}\,\to \,0
\,\quad\mbox{as }\,n\to\infty,
\end{align}
and the validity of \eqref{adjoint} is shown. Next, we take the limit as $n\to\infty$ in the variational
inequality \eqref{vugal}, written for $\alpha=\alpha_n$. Employing \eqref{conny1}, \eqref{conny2}, and
\eqref{conny5}, we readily see that \eqref{vug} is fulfilled. This concludes the proof of the assertion.\QED

\vspace{3mm}
\Brem
Unfortunately, we are unable to derive a complementarity slackness condition for the adjoint variables.
Indeed, although we have the inequality 
\begin{align}
&\langle(\varphi(\alpha_n)h''(\rho^{\alpha_n})q^{\alpha_n},\varphi(\alpha_n)h''(\rho_\Gamma^{\alpha_n})q_\Gamma^{\alpha_n}),
(q^{\alpha_n},q_\Gamma^{\alpha_n})\rangle_{{\cal Z}_0}\non\\[2mm]
&=\,\intQ \varphi(\alpha_n)h''(\rho^{\alpha_n})\left|q^{\alpha_n}\right|^2\,+\intS
\varphi(\alpha_n)h''(\rho_\Gamma^{\alpha_n})\left|q_\Gamma^{\alpha_n}\right|^2\,\ge\,0\quad\forall\,n\in\enne,
\end{align}
the convergence properties \eqref{conny4} and \eqref{conny6} are not strong enough to guarantee that
$\langle (\lambda,\lambda_\Gamma),(q,\qG)\rangle_{{\cal Z}_0}\,\ge\,0$.
\Erem

\vspace{3mm}
\Brem
Obviously, the adjoint variables are not uniquely determined. It thus may well happen that for different sequences $\alpha_n\searrow0$
different limits are approached. However, \pier{the weak-star limit $\eta$ in \eqref{conny5} must satisfy the variational inequality \eqref{vug}.}
\Erem

\vspace{3truemm}


\Begin{thebibliography}{10}

\bibitem{Barbu}
V. Barbu\pier{:}
{\em Necessary conditions for nonconvex distributed control problems governed 
by elliptic variational inequalities}.
J. Math. Anal. Appl. {\bf 80} (1981), 566-597.

\bibitem{Brezis}
H. Brezis\pier{:
``Op\'erateurs maximaux monotones et semi-groupes de contractions
dans les espaces de Hilbert''}.
North-Holland Math. Stud. {\bf 5},
North-Holland,
Amsterdam,
1973.

\bibitem{CFGS1}
P. Colli, \pier{M.H. Farshbaf-Shaker, G. Gilardi,} J. Sprekels: 
{\em Optimal boundary control of a 
viscous Cahn-Hilliard system with dynamic boundary condition 
and double obstacle potentials}.
SIAM J. Control Optim. {\bf 53} (2015), 2696-2721.

\bibitem{CFGS2}
P. Colli, M.H. Farshbaf-Shaker, G. Gilardi, J. Sprekels: {\em Second-order 
analysis of a boundary control problem for the viscous Cahn-Hilliard 
equation with dynamic boundary conditions.}
\pier{Ann. Acad. Rom. Sci. Ser. Math. Appl.} {\bf 7} (2015), 41-66.

\bibitem{CFS}
P. Colli, M.H. Farshbaf-Shaker, J. Sprekels: {\em A deep quench approach 
to the optimal control of an Allen--Cahn equation with dynamic boundary 
conditions and double obstacles.}
Appl. Math. Optim. {\bf 71} (2015), 1-24.  

\pier{\bibitem{CF2} 
P.\ {C}olli, T.\ {F}ukao: 
{\em Equation and dynamic boundary condition of 
Cahn--Hilliard type with singular potentials.}
Nonlinear Anal. {\bf 127} (2015), 413-433.%
}


\bibitem{CGRS1}
P. Colli, G. Gilardi, E. Rocca, J. Sprekels:
{\em Optimal distributed control of a diffuse interface model of tumor growth.} 
Nonlinearity \pier{{\bf 30} (2017), 2518-2546.}

\bibitem{CGS1}
P. Colli, G. Gilardi, J. Sprekels: {\em A boundary control problem 
for the pure Cahn--Hilliard equation 
with dynamic boundary conditions.}
Adv. Nonlinear Anal. {\bf 4} (2015), 311-325.

\bibitem{CGS2}
P. Colli, G. Gilardi, J. Sprekels: {\em Distributed optimal control 
of a nonstandard nonlocal phase field system.}
AIMS \pier{Mathematics} {\bf 1} (2016), 246-281.

\bibitem{CGS3}
P. Colli, G. Gilardi, J. Sprekels: {\em A boundary control problem 
for the viscous Cahn--Hilliard equation with 
dynamic boundary conditions.}
Appl. Math. Optim. {\bf 72} (2016), 195-225.

\pier{\bibitem{CGS3bis}
P. Colli, G. Gilardi, J. Sprekels: {\em Recent results 
on the Cahn--Hilliard equation 
with dynamic boundary conditions.} Vestn. Yuzhno-Ural. 
Gos. Univ., Ser. Mat. Model. Program. {\bf 10} (2017), 5-21.}

\bibitem{CGS4}
P. Colli, G. Gilardi, J. Sprekels: {\em Distributed optimal 
control of a nonstandard nonlocal phase field 
system with double obstacle potential.}
Evol. Equ. Control Theory {\bf 6} (2017), 35-58.

\bibitem{CGS13}
P. Colli, G. Gilardi, J. Sprekels: {\em On a Cahn--Hilliard system with convection
and dynamic boundary conditions}. \pier{Preprint arXiv:1704.05337 [math.AP] (2017), 
pp.~1-34.}

\bibitem{CGS14}
P. Colli, G. Gilardi, J. Sprekels: {\em Optimal velocity control 
of a viscous Cahn--Hilliard system with convection and dynamic boundary conditions}. \pier{Preprint arXiv:170X.????? [math.AP] (2017), pp.~1-29.}

\bibitem{CS10}
P. Colli, J. Sprekels: Optimal boundary control of a nonstandard 
Cahn--Hilliard system with dynamic 
boundary condition and double obstacle inclusions. 
\pier{To appear in ``Solvability, Regularity, Optimal Control 
of Boundary Value Problems for PDEs'', 
P.~Colli, A.~Favini, E.~Rocca, G.~Schimperna, J.~Sprekels~(eds.), 
Springer INdAM Series (see also the WIAS Preprint No. 2370, 
Berlin 2017, pp.~1-29).} 
	
\pier{\bibitem{DZ}
N. Duan, X. Zhao: {\em Optimal control for the multi-dimensional 
viscous Cahn--Hilliard equation.}
Electron. J. Differential Equations 2015, No. 165, 13 pp.}


\bibitem{FRS}
S. Frigeri, E. Rocca, J. Sprekels: {\em Optimal distributed control 
of a nonlocal Cahn--Hilliard/Navier--Stokes 
system in two dimensions.}
SIAM J. Control. Optim. {\bf 54} (2016), 221-250. 

\pier{%
\bibitem{FY} 
T. Fukao, N. Yamazaki:
{\em A boundary control problem for the equation 
and dynamic boundary condition of Cahn--Hilliard type.} 
To appear in ``Solvability, Regularity, Optimal Control 
of Boundary Value Problems for PDEs'', 
P.~Colli, A.~Favini, E.~Rocca, G.~Schimperna, J.~Sprekels~(eds.), 
Springer INdAM Series.}

\bibitem{GMS}
G. Gilardi, A. Miranville, G. Schimperna: {\em On the Cahn--Hilliard 
equation with irregular potentials
and dynamic boundary conditions.} Commun. Pure Appl. Anal. {\bf 8} (2009), 881-912.

\bibitem{HHKW}
M. Hinterm\"uller, M. Hinze, C. Kahle, T. Kiel: {\em A goal-oriented 
dual-weighted adaptive finite element approach for the optimal 
control of a nonsmooth Cahn--Hilliard--Navier--Stokes system.}
WIAS Preprint No. 2311, Berlin 2016, \pier{pp.~1-27.}

\bibitem{HKW}
M. Hinterm\"uller, T. Kiel, D. Wegner: {\em Optimal control of a semidiscrete 
Cahn--Hilliard--Navier--Stokes system with non-matched fluid densities}. 
\pier{SIAM J. Control Optim. {\bf 55} (2017), 1954-1989.}

\bibitem{hw}
M. Hinterm\"uller, D. Wegner: {\em
Distributed optimal control of the Cahn--Hilliard system 
including the case of a double-obstacle homogeneous free energy density}. 
SIAM J. Control Optim. {\bf 50} (2012), 388-418. 

\bibitem{HW1}
M. Hinterm\"uller, D. Wegner: {\em Optimal control of a semidiscrete 
Cahn--Hilliard--Navier--Stokes system}.
SIAM J. Control Optim. {\bf 52} (2014), 747-772.

\bibitem{HW2}  
M.  Hinterm\"uller,  D.  Wegner: {\em Distributed  and  boundary  control  problems  
for  the  semidiscrete Cahn--Hilliard/Navier--Stokes system with nonsmooth 
Ginzburg--Landau energies}. 
Isaac Newton Institute Preprint Series No. NI14042-FRB (2014), \pier{pp.~1-29.}

\bibitem{Kudla}
C. Kudla, A.T. Blumenau, F. B\"ullesfeld, N. Dropka, C. Frank-Rotsch, 
F. Kiessling, O. Klein, P. Lange, W. Miller,
U. Rehse, U. Sahr, M. Schellhorn, G. Weidemann, M. Ziem, 
G. Bethin, R. Fornari, M. M\"uller, J. Sprekels,
V. Trautmann, P. Rudolph: 
{\em Crystallization of 640 kg mc-silicon ingots under traveling magnetic
field by using a heater-magnet module.} J. Crystal Growth {\bf 365} (2013), 54-58.


\bibitem{RS}
E. Rocca, J. Sprekels: {\em Optimal distributed control of a nonlocal 
convective Cahn--Hilliard equation by the velocity in three dimensions}. 
SIAM J. Control Optim. {\bf 53} (2015), 1654-1680.

\bibitem{Simon}
J. Simon: {\em Compact sets in the space $L^p(0,T; B)$}.
Ann. Mat. Pura Appl. \pier{(4)}
{\bf 146} (1987), 65-96.

\pier{\bibitem{TM}
T. Tachim Medjo: {\em Optimal control of a Cahn--Hilliard--Navier--Stokes 
model with state constraints.}
 J. Convex Anal.  {\bf 22}  (2015),  1135-1172.}

\bibitem{WaNa} Q.-F. Wang, S.-i. Nakagiri: 
{\em Weak solutions of Cahn--Hilliard equations 
having forcing terms and optimal control problems}. 
Mathematical models in functional equations (Japanese) (Kyoto, 1999), 
S\={u}rikaisekikenky\={u}sho K\={o}ky\={u}roku No. 1128 (2000), 172--180.

\bibitem{ZL1}
X.P. Zhao, C.C. Liu: {\em Optimal control of the convective 
Cahn--Hilliard equation}. Appl. Anal. {\bf 92} (2013), 1028-1045.

\bibitem{ZL2} 
X.P. Zhao, C.C. Liu: {\em Optimal control of the convective 
Cahn--Hilliard equation in 2D case}. Appl. Math. Optim. {\bf 70} (2014), 61-82.
		
\pier{\bibitem{Zh}
J. Zheng: {\em Time optimal controls of the Cahn--Hilliard equation with internal control.}  Optimal Control Appl. Methods  {\bf 36}  (2015), 566-582.}

\bibitem{ZW}
J. Zheng, Y. Wang: {\em Optimal control problem for Cahn--Hilliard 
equations with state constraint}.
J. Dyn. Control Syst. {\bf 21} (2015), 257-272.

\End{thebibliography}


\End{document}
